\documentclass[a4paper,11pt,fleqn]{article}
\usepackage[utf8]{inputenc}
\usepackage[T1]{fontenc}    
\usepackage[ruled,vlined]{algorithm2e}
\usepackage{amsmath,amsfonts,amssymb,amsthm}
\usepackage{bm,bbm}
\usepackage{mathtools}
\usepackage{mathrsfs} 
\usepackage{thm-restate,thmtools}
\usepackage{soul}
\usepackage{algorithmic}
\usepackage{nicefrac}
\usepackage{color}
\usepackage{xcolor,colortbl}
\usepackage{nicefrac}  
\usepackage{booktabs}
\usepackage{adjustbox}
\usepackage{hyperref}
\usepackage{url}
\usepackage{cleveref}
\usepackage[flushleft]{threeparttable}
\usepackage{multirow}
\usepackage[shortlabels]{enumitem}
\usepackage{pifont}
\usepackage{microtype}      
\usepackage{natbib}
\usepackage{doi}
\usepackage{authblk}
\usepackage{tikz}
\usepackage{tablefootnote}
\usepackage{needspace}
\usepackage{euscript}
\usepackage{charter}
\usepackage{array}
\usepackage[normalem]{ulem}
\usepackage{lipsum}
\usepackage{graphicx}
\graphicspath{ {./images/} }

\newtheorem{theorem}{Theorem}[section]

\newtheorem{lemma}[theorem]{Lemma}

\newtheorem{definition}[theorem]{Definition}
\newtheorem{assumption}{Assumption}

\newcommand{\cA}{\mathcal{A}}
\newcommand{\cH}{\mathcal{H}}

\newcommand{\cS}{\mathcal{S}}
\newcommand{\cX}{\mathcal{X}}

\newcommand{\fo}{\mathfrak{o}}

\newcommand{\bbE}{\mathbb{E}}
\newcommand{\bbN}{\mathbb{N}}
\newcommand{\bbR}{\mathbb{R}}
\newcommand{\bbP}{\mathbb{P}}

\newcommand{\1}{\mathbf{1}}
\newcommand{\cF}{\mathcal{F}}

\newcommand{\prox}{\operatorname{prox}}

\newcommand{\as}{\mathrm{a.s.}}

\newcommand{\norm}[1]{\left\lVert #1\right\rVert}

\newcommand{\ip}[2]{\left\langle #1,#2\right\rangle}
\newcommand{\proj}{\operatorname{\Pi}}
\newcommand{\sgn}{\operatorname{sgn}}

\definecolor{labelkey}{rgb}{0,0.08,0.45}
\definecolor{refkey}{rgb}{0,0.6,0.0}
\definecolor{Brown}{rgb}{0.45,0.0,0.05}
\definecolor{dgreen}{rgb}{0.00,0.49,0.00}
\definecolor{dblue}{rgb}{0,0.08,0.75}

\hypersetup{
colorlinks,
linkcolor=dgreen,
citecolor=dblue
}

\numberwithin{equation}{section}

\title{ { \sffamily Nearly Tight Almost Sure Convergence Rates for the Last iterate of SGD} } 

\author[1]{Weiming Ou}
\author[2]{Andrea Paudice}
\author[1, 3]{Xiao Wang}

\affil[1]{\footnotesize School of Computing and Artificial Intelligence, Shanghai University of Finance and Economics, No.777 Guoding Road, 200433 Shanghai, China}
\affil[2]{\footnotesize Department of Computer Science, Aarhus University, Åbogade 34, 8200 Aarhus, Denmark}
\affil[3]{\footnotesize MoE Key Laboratory of Interdisciplinary Research of Computation and Economics, Shanghai University of Finance and Economics, No.777 Guoding Road, 200433 Shanghai, China}

\date{}

\begin{document}
\maketitle

\begin{abstract}
We study the almost sure convergence of the last-iterate optimization error of \emph{Stochastic (sub)Gradient Descent} (SGD) for convex objectives. For smooth objectives, we establish convergence rates of order $k^{-1/3}$ when a stationary minimizer exists, and $k^{-1/4}$ otherwise, up to polylogarithmic factors, where $k$ denotes the iteration counter. We further show that the latter rate extends to nonsmooth convex objectives. All these rates are achieved by deterministic, \emph{anytime} stepsize policies. Finally, we establish matching lower bounds showing that the exponents $-1/3$ and $-1/4$ are optimal, up to polylogarithmic factors, for the corresponding problem classes. 
\end{abstract}

\vspace{1ex}
\allowdisplaybreaks

\section{Introduction}
\emph{Stochastic (sub)Gradient Descent} (SGD) \citep{Robbins1951} and its nonsmooth generalization~\citep{Ermoliev1969} are the workhorse of modern machine learning: from linear prediction to deep neural networks, large-scale models are trained by cheap noisy gradient steps on samples or mini-batches \citep{Bottou2010,LeCun2015,Bottou2018}. Its classical guarantees, however, concern \emph{averaged} iterates, which attain optimal rates by suppressing the fluctuations of individual iterates \citep{Polyak1992,Nemirovski2009}. In practice, it is the \emph{last iterate} that is monitored, deployed, and used to decide when to stop. This mismatch has shifted the community's attention from averaged to last-iterate guarantees, a question raised explicitly by \citet{Shamir2012} and since addressed by a growing literature on sharp last-iterate rates in expectation and with high probability \citep{Moulines2011,Shamir2013,Harvey2019,Jain2021,Liu2024b}.

\noindent \noindent Nonasymptotic analyses provide expectation or high-probability bounds at prescribed iterations, as well as \emph{time-uniform} high-probability bounds that support repeated monitoring and data-dependent stopping \citep{Chen2026,Li2026}. Almost-sure rates address a complementary question: \emph{how fast does the optimization error decay along almost every realization of a single run?} Such a rate holds at every sufficiently large iteration on a probability-one event, with a path-dependent constant and starting time, but need not provide an explicit finite-time confidence bound. Pointwise high-probability bounds for a single horizon-free run can also yield almost-sure rates through the Borel--Cantelli lemma. When their dependence on $1/\delta$ is polylogarithmic (for bounds holding with probability at least $1-\delta$), choosing summable failure probabilities costs only polylogarithmic factors in the iteration counter, as in light-tailed analyses \citep{Harvey2019,Eldowa2024}. Under second-moment assumptions alone, however, polynomial dependence on $1/\delta$ can arise \citep{Sadiev2023}, and this direct conversion may worsen the convergence exponent or yield no decaying bound. This motivates a direct study of optimal almost-sure last-iterate rates under such assumptions: do rare, late stochastic fluctuations impose a fundamental limitation on convergence?


\noindent Even in the basic convex setting, almost sure convergence analyses are scarce, and the optimal rate remains unknown. Beyond classical 
\emph{qualitative} (i.e., without an explicit rate) convergence results \citep{Robbins1951,Robbins1971}, \citet{Sebbouh2021} established almost-sure rates for weighted averages, while \citet{Liu2024a} proved last-iterate rates $\mathcal{O}(k^{-1/3+\beta})$ with any $\beta\in (0,1/3)$ for smooth convex objectives under a standard variance assumption on the gradient noise \citep{Khaled2020}, with further results in related settings \citep{Hudiani2025,Weissmann2024,Chen2026}. For nonsmooth convex objectives, qualitative almost-sure convergence of the stochastic subgradient method is classical \citep{Ermoliev1969,Polyak1987,KushnerYin2003}, and almost-sure rates for weighted averages follow from the Robbins--Siegmund theorem \citep{Robbins1971}, by the argument of \citet{Sebbouh2021}, which uses smoothness only to bound the second moment of the oracle. Notwithstanding these advances, there are still fundamental gaps unresolved. In particular, for smooth objectives, it remains unknown whether the exponent $-1/3$ is intrinsic and how it changes if we add constraints. For nonsmooth objectives, to the best of our knowledge, no almost-sure rate is known for its \emph{last iterate} under a uniform bound on the second moment of the stochastic subgradients.

\paragraph{Contributions.}
We establish nearly matching lower and upper bounds for the almost-sure convergence of the last iterate of SGD only requiring conditions on the moment of the stochastic (sub)gradients.
\begin{itemize}[leftmargin=*]
\item \textbf{Lower bounds.} For a broad class of stepsize policies, including all deterministic (bounded) schedules and certain randomized and adaptive policies, we prove almost-sure last-iterate lower bounds of order $k^{-1/3}$ for smooth convex objectives with a stationary minimizer (\Cref{thm:smooth_lb}~\ref{itm:smooth_lb_fast}), and $k^{-1/4}$ for smooth convex objectives without a stationary minimizer and for nonsmooth convex objectives, up to logarithmic factors (\Cref{thm:smooth_lb}~\ref{itm:smooth_lb_slow} and \Cref{thm:nonsmooth_lb}). These bounds show that the exponent $-1/3$ in existing smooth analyses is intrinsic and that, without a stationary minimizer, smoothness does not improve the worst-case exponent over the nonsmooth case.
\item \textbf{Matching upper bounds.} We give deterministic, \emph{anytime} stepsize policies attaining these rates up to polylogarithmic factors (\Cref{thm:smooth_ub,thm:nonsmooth_ub}). For smooth objectives with a stationary minimizer, this sharpens the $k^{-1/3+\beta}$ rates of \citet{Liu2024a}, valid for arbitrarily small $\beta>0$, by replacing the polynomial slack with polylogarithmic factors. For nonsmooth convex objectives, we provide, to the best of our knowledge, the first almost-sure last-iterate rate under a bounded conditional second-moment assumption alone.
\end{itemize}

\subsection{Related Works.}
\paragraph{Expected Error Analysis.}
For smooth convex objectives, \citet{Moulines2011} obtained a {last-iterate} rate of order $k^{-1/3}$ with stepsizes proportional to $k^{-2/3}$, under suitable second-moment assumptions on the stochastic gradients. Under uniformly bounded gradient-noise variance, \citet{Liu2024b} established a {last-iterate} rate of order $k^{-1/2}$, up to logarithmic factors, within a framework also covering nonsmooth objectives and general domains. This nearly matches the minimax rate under bounded-variance noise \citep{Nemirovsky1983,Agarwal2009}. For nonsmooth convex objectives on bounded domains, with uniformly bounded stochastic subgradient second moments, \citet{Shamir2013} proved a {last-iterate} rate of order $(\log k) \cdot k^{-1/2}$ with horizon-independent
stepsizes proportional to $k^{-1/2}$.

\paragraph{Error Tail Bounds.}
Markov's inequality converts expected error bounds into high-probability guarantees with a $1/\delta$ overhead. Under almost-surely bounded oracle noise, \citet{Harvey2019} obtained a last-iterate bound of order $(\log k) \cdot k^{-(1/2)} \cdot \log(1/\delta)$ for nonsmooth convex Lipschitz objectives on bounded domains, with extensions to sub-Gaussian noise. For both smooth and nonsmooth convex objectives on general domains, \citet{Liu2024b} established last-iterate rates of order $k^{-1/2}$, up to logarithmic factors in $k$ and $1/\delta$, under conditionally sub-Gaussian noise, including with horizon-independent stepsizes. Their extended version also covers conditionally sub-Weibull noise~\citep{Vladimirova2020}, with polylogarithmic dependence on $1/\delta$. Similar results have also been obtained by \cite{Eldowa2024}, who also consider, at a cost of a worse depedence on $1/\delta$, the case of polynomially tailed noise. Under bounded variance alone, however, \citet{Sadiev2023} showed that constant-stepsize SGD can require polynomial dependence on $1/\delta$, even for smooth strongly convex objectives.

\section{Problem Setting}
In this work, we consider the following problem:
\begin{equation}\label{eq:problem}
\underset{x \in \cX}{\operatorname{minimize}}\, F(x)  \tag{P}
\end{equation}
with access only to a stochastic first-order oracle, whose response at step $k$ is denoted by $\widehat{\nabla}_k F(X_k)$.

\noindent To make \eqref{eq:problem} tractable, we make the following standard assumptions. 
\begin{assumption}[Convexity]\label{ass:convexity}
$\cH$ is a real separable Hilbert space and $F\colon \cH \rightarrow \bbR$ is convex. 
\end{assumption}
\begin{assumption}[Feasibility]\label{ass:domain}
$\cX \subseteq \cH$ is closed and convex, and $\mathcal{X}_{\star}\coloneqq \arg \min_{x \in \cX} F(x) \neq \varnothing$. We define $F_{\star} \coloneqq  \min_{x \in \cX} F(x)$.
\end{assumption}
\paragraph{Notation.}
Throughout we mainly use uppercase letters for random iterates; their coordinates and oracle noises may use lowercase letters. For a set $\cA$, define the extended-valued indicator $\iota_{\cA}(a)=0$ when $a \in \cA$ and $+\infty$ otherwise. We use $k$ to denote the step counter, and therefore it will always be a natural number. All random variables are defined on a common probability space $(\Omega, \cF, \bbP)$ equipped with a filtration $(\cF_k)_{k \geq 1}$. Here, $\cF_k$ represents the information available before the oracle call at iteration $k$, including any algorithmic randomness used to select the step-size. The iterate $X_k$ is $\cF_k$-measurable, and the oracle response $\widehat{\nabla}_k F(X_k)$ is $\cF_{k+1}$-measurable. The conditional distribution of this response may depend on the information in $\cF_k$. We define the conditional expectation by $\bbE_k[\cdot] \coloneqq \bbE[\cdot | \cF_k]$\footnote{The conditional expectation of Hilbert valued random variables is meant in the Bochner sense \citep{HytonenEtAl2016}.} and the conditional probability by $\bbP_k(\cdot) := \bbP(\cdot|\cF_k)$ for $k \geq 1$. {Since conditional expectation and conditional probability are defined up to null sets, all formulas involving $\bbE_k$ and $\bbP_k$ are understood almost surely.} Denote by $\partial f$ the subdifferential of any proper convex $f$. With abuse of notation, denote by $\widehat{\nabla}_k F(X_k)$ the stochastic gradient (resp. subgradient) oracle in the smooth (resp. nonsmooth) setting. Given a sequence of random variables $(Y_k)_{k \geq 1}$ and a sequence of reals $(y_k)_{k \geq 1}$, when we write $Y_k \equiv y_k$ we mean that, for every $k \geq 1$, $Y_k(\omega) = y_k$ for every $\omega \in \Omega$. For $g>0$, we write $f(n)=o(g(n))$ if $f(n)/g(n)\to0$ as $n\to\infty$, understood pathwise for random sequences. Finally, for $x \in \cX$, let us define the optimization error as $\Delta(x) \coloneqq F(x) - F_{\star}$.

\paragraph{Stochastic (Sub)Gradient Descent (SGD).} Let $x_1 \in \mathcal{X}$ and $(\Gamma_k)_{k \geq 1}$ be a (possibly random) sequence of nonnegative, finite step-sizes. We consider the following stochastic recurrence
\begin{align}
\label{eq:sgd}
\forall k \geq 1, \quad X_{k+1} =\proj_\mathcal X\left( X_k - \Gamma_k \cdot \widehat{\nabla}_k F(X_k) \right),
\end{align}
where $X_1 = x_1$ almost surely and $\proj_{\mathcal X}$ is the orthogonal projection operator onto $\mathcal X$. Specifically, we consider \emph{predictable} step-size policies as defined by the following.  
\begin{definition}\label{def:predictable_ss}
Let $\alpha > 0$. A sequence of positive real-valued random variables $(\Gamma_k)$ is $\alpha$-bounded and predictable w.r.t. the filtration $(\cF_k)_{k \geq 1}$ if for every $k \geq 1$, $\Gamma_k$ is $\cF_k$-measurable and $\Gamma_k \in (0, \alpha]$ almost surely. We denote the collection of all such sequences with $\cS_{\alpha}$. We also consider the class $\cS_{\mathrm{fin}} \coloneqq \{(\Gamma_k)_{k \geq 1}: \text{ for every } k \geq 1, \Gamma_k \text{ is $\cF_k$-measurable and } 0 \leq \Gamma_k < +\infty\}$, and observe that $\cS_{\alpha} \subseteq \cS_{\mathrm{fin}}$ for every $\alpha > 0$.
\end{definition}
The definition of $\cS_{\alpha}$ includes all deterministic schedules, given by $\Gamma_k \equiv \gamma_k$ for $(\gamma_k)_{k\ge1}$ satisfying
$0<\gamma_k\le\alpha$ for every $k\ge1$. In particular, it includes positive nonincreasing schedules with $\gamma_1\le\alpha$. It also includes bounded adaptive rules depending on past oracle responses, such as the scalar adaptive schedule based on past gradients in \cite{Li2019}, as well as randomized policies whose internal randomization is available before
the current oracle call. In contrast, gradient clipping \citep{Gorbunov2020} and the usual AdaGrad-Norm \citep{Ward2020} update are generally not covered: when expressed in the form~\eqref{eq:sgd} with the original oracle response, their effective scalar step-sizes also depend on the current response. See Appendix~\ref{app:details_predictable_ss} for details. 

\noindent Our focus is on the convergence rate of $\Delta(X_k)$ as captured by the following definition.
%
%

\begin{definition}[Almost sure convergence]\label{def:as_convergence}
Let $(\rho_k)_{k \geq 1}$ be a sequence of positive real numbers such that $\rho_k \to 0$ for $k \to +\infty$, and $(Y_k)_{k \geq 1}$ a sequence of real valued random variables. We say that $(Y_k)_{k \geq 1}$ converges to $0$ almost surely with rate $\rho_k$ and write $Y_k \in \fo_{a.s.}(\rho_k)$ iff 
\begin{align*}
\bbP \bigg( \underset{k \to \infty}{\lim} \frac{Y_k}{\rho_k} = 0 \bigg) = 1. 
\end{align*}
If this does not hold, then we write $Y_k \notin \fo_{a.s.}(\rho_k)$.

\end{definition}
For consistency with our definition of almost-sure convergence rates, we state the lower bounds in Theorem~\ref{thm:smooth_lb} and~\ref{thm:nonsmooth_lb} in terms of $\fo_{a.s.}$ notation. However, our constructions establish a stronger conclusion. Specifically, for the corresponding hard instances and associated rate sequences $(\rho_k)_{k\ge 1}$, we show that 
\begin{align*}
\bbP \bigg(\underset{k \to \infty}{\limsup} \frac{\Delta(X_k)}{\rho_k} = \infty \bigg) = 1,
\end{align*}
which clearly implies $\Delta(X_k) \notin \fo_{a.s.}(\rho_k)$. {In the proofs, we write $r_k=1/\rho_k$ for the reciprocal of the rate under consideration, so that $\Delta(X_k)\in\fo_{a.s.}(\rho_k)$ if and only if $r_k \cdot \Delta(X_k)\to0$ almost surely, while the lower bounds establish $\limsup_k r_k \cdot \Delta(X_k)=\infty$ almost surely.}

\section{Smooth Optimization}\label{sec:smooth_obj}
In this section, we present our results on smooth objectives. We start from adopting the following standard notion of smoothness.

\begin{assumption}[Smoothness]\label{ass:smoothness}
$F\colon \cH \rightarrow \bbR$ is differentiable and there exists $L > 0$ such that for every $x, y\in \cH$, it holds that \(\|\nabla F(x)-\nabla F(y)\| \leq L \cdot \|x-y\|\).
%
\end{assumption}
Next, we specify the statistical model of the first-order oracle. We use a conditional version of the $(a, b, c)$ condition introduced by \citet{Khaled2020}.
\begin{assumption}\label{ass:abc_condition}
There exist $a, b, c \ge 0$ such that for each $k \geq 1$, $\widehat{\nabla}_k F(X_k)$ satisfies the conditions: 
\begin{enumerate}[label=\textnormal{(\roman*)}]
\item (Unbiasedness) $\bbE_k[\widehat{\nabla}_k F(X_k)] = \nabla F(X_k)$.
\item\label{itm:abc-moment} (abc condition) $\mathbb{E}_k[\|\widehat{\nabla}_k F(X_k)\|^2]\leq a \cdot \Delta(X_k) + b \cdot \|\nabla F(X_k)\|^2 + c$.
\end{enumerate}
\end{assumption}
Finally, we introduce the following assumption as a relaxation of the unconstrained version of \eqref{eq:problem}.  
\begin{assumption}\label{ass:stationary_minimizer}
There exists $x_{\star} \in \cX_\star$ such that $\nabla F(x_{\star})=0.$
\end{assumption}
The above is always satisfied in the case of unconstrained optimization (i.e. $\cX = \cH$) or when there exists at least one minimizer in the interior of the feasible set. We can now present the main results of this section, starting from the following lower bounds.
\begin{theorem}\label{thm:smooth_lb}
Let $\alpha > 0$, and $L > 0$.  
\begin{enumerate}[label=\textnormal{(\roman*)}]
\item\label{itm:smooth_lb_slow} Let $\cH = \bbR^2$ and $\cX = [0,1] \times [-2,2]$. Then, there exists an objective $F\colon \cH \rightarrow \bbR$ satisfying Assumptions~\ref{ass:convexity}-\ref{ass:smoothness} with the prescribed $L$, a first order oracle satisfying Assumption~\ref{ass:abc_condition} and $x_1 \in \cX$ such that for every $(\Gamma_k)_{k \geq 1} \in \cS_{\mathrm{fin}}$ and every $p > 0$, the iterates generated by eq.~\eqref{eq:sgd} satisfy $\Delta(X_k) \notin \fo_{a.s.}\left((\log(e+k))^{-p} \cdot k^{- 1/4 }\right)$.

\item\label{itm:smooth_lb_fast} Let $\cH = \cX = \bbR^2$. Then, there exists an objective $F\colon \cH \rightarrow \bbR$ satisfying Assumptions~\ref{ass:convexity}-\ref{ass:smoothness}, and \ref{ass:stationary_minimizer} on $\cX$ with the prescribed $L$, a first order oracle satisfying Assumption~\ref{ass:abc_condition} and $x_1 \in \cX$ such that for every $(\Gamma_k)_{k \geq 1} \in \cS_{\alpha}$ and every $p > 0$, the iterates generated by eq.~\eqref{eq:sgd}  satisfy $\Delta(X_k) \notin \fo_{a.s.}((\log(e+k))^{-p} \cdot k^{- 1/3})$.
\end{enumerate}
\end{theorem}

\noindent For unconstrained problems over $\bbR^d$, the analysis of \citet{Liu2024a} implies that $\Delta(X_k)\in\fo_{a.s.}(k^{-1/3+\beta})$ for every $\beta\in(0,1/3)$, using a suitable deterministic step-size policy. Since the construction in \Cref{thm:smooth_lb}~\ref{itm:smooth_lb_fast} already applies with $\cX=\cH=\bbR^2$, together these results identify $-1/3$ as the optimal exponent for this class of problems. This yields a polynomial separation between the last iterate and a suitable average iterate. Indeed, the latter converges to $0$ with a rate of $k^{-1/2+\beta}$ for every $\beta\in(0,1/2)$ under standard second-moment assumptions \citep{Sebbouh2021}. Thus, averaging can improve the achievable polynomial exponent, rather than merely yielding a logarithmic-factor improvement as in existing non-asymptotic analyses. There is also a separation from the non-asymptotic behavior of the last iterate itself. Under uniformly bounded noise variance, \citet{Liu2024b} establish an expected error of order $(\log k) \cdot k^{-1/2}$. Our lower bound shows that such expectation guarantees cannot, in general, be matched in the almost-sure sense, even by allowing adaptive (but predictable) step-sizes. The distinction is therefore that controlling the expected error at each time and controlling its decay along an entire sample path lead to different polynomial barriers. \Cref{thm:smooth_lb}~\ref{itm:smooth_lb_slow} reveals an additional role of the geometry of the constrained optimum. Without Assumption~\ref{ass:stationary_minimizer}, even globally smooth convex objectives can preclude every almost-sure last-iterate rate with polynomial exponent beating $-1/4$. Hence, smoothness alone does not preserve the $-1/3$ exponent when the feasible set contains no stationary minimizer.

\noindent The next result provides upper bounds for both the settings in the items of \Cref{thm:smooth_lb}. 
%
%
\begin{theorem}\label{thm:smooth_ub}
Fix $\eta > 0$ and suppose that Assumptions \ref{ass:convexity}-\ref{ass:abc_condition} hold.
\begin{enumerate}[label=\textnormal{(\roman*)}]
\item\label{itm:smooth_ub_slow} There exist $\alpha, t_0 > 0$ and an $\alpha$-bounded predictable step-size sequence $(\Gamma_k)_{k \geq 1}$ such that the iterates generated by eq.~\eqref{eq:sgd} satisfy $\Delta(X_k) \in \fo_{a.s.}([\log(k+t_0)]^{1/2+\eta} \cdot (k+t_0)^{-1/4})$.

\item\label{itm:smooth_ub_fast} Suppose that Assumption \ref{ass:stationary_minimizer} also holds. There exist $\alpha, t_0 > 0$ and an $\alpha$-bounded predictable step-size sequence $(\Gamma_k)_{k \geq 1}$, depending on $t_0$, such that the iterates generated by eq.~\eqref{eq:sgd} satisfy $\Delta(X_k) \in \fo_{a.s.}([\log(k+t_0)]^{1/3+\eta} \cdot (k+t_0)^{-1/3})$.
\end{enumerate}
\end{theorem}
\Cref{thm:smooth_ub}~\ref{itm:smooth_ub_slow} establishes, to the best of our knowledge, the first almost sure last-iterate guarantee for general smooth convex constrained problems under \Cref{ass:abc_condition}, without requiring a stationary minimizer. In particular, the resulting rate has exponent $-1/4$ and explicit logarithmic correction $[\log(k+t_0)]^{1/2+\eta}$. Along with \Cref{thm:smooth_lb}~\ref{itm:smooth_lb_slow}, it identifies $-1/4$, up to poly-logarithmic factors, as the optimal exponent for the problem class. \Cref{thm:smooth_ub}~\ref{itm:smooth_ub_fast} sharpens the guarantee of \citet{Liu2024a} in the unconstrained setting, replacing an arbitrarily small polynomial loss with the explicit logarithmic correction $[\log(k+t_0)]^{1/3+\eta}$. Moreover, the result extends to projected SGD whenever the feasible set contains a stationary minimizer. Thus, the $-1/3$ exponent is not restricted to unconstrained optimization. 
An inspection of the proofs also reveals that the adopted step-sizes policies are deterministic and horizon-independent, whereas adaptive predictable scalar step-sizes cannot improve the exponents in the worst case. The distinction between the two regimes therefore reflects a genuine role of the constrained optimum, rather than merely a limitation of the upper-bound analysis. On the other hand, the optimality of the logarithmic corrections remains unresolved by our bounds.

\subsection{Analysis}\label{sec:smooth-sketch}
Below we provide the proof of Theorem~\ref{thm:smooth_lb}~\ref{itm:smooth_lb_slow}, while deferring the proof of Theorem~\ref{thm:smooth_lb}~\ref{itm:smooth_lb_fast} and \Cref{thm:smooth_ub} to Appendices ~\ref{app:smooth-endpoint} and ~\ref{app:upper-smooth}, respectively. All the referred lemmas can be found in Appendices ~\ref{app:common tools} and ~\ref{app:smooth-lb} along with their proofs.

\paragraph{Overview.} The two lower bounds of \Cref{thm:smooth_lb} arise from the same conflict: the stepsizes must be large enough to make progress, but small enough to control fresh noise. We place these requirements in separate coordinates of a two-dimensional objective. One coordinate contains increasingly flat regions, so reaching smaller objective values requires sufficient cumulative movement. The other is perturbed by independent finite-variance innovations, and the same scalar stepsize controls both updates. The geometry of this noisy coordinate distinguishes the two cases. At a nonstationary boundary minimizer, an inward displacement of size $h$ produces a linear objective gap; at a stationary minimizer, we use a quadratic coordinate, so a noise displacement of size $h$ forces an error of order $h^2$ at the current or next iterate. This difference underlies the barriers $k^{-1/4}$ and $k^{-1/3}$. This argument identifies the critical exponents, but obtaining almost-sure lower bounds requires the obstruction to recur infinitely. Repeating equally strong shocks at infinitely many scales would violate the variance constraint in \Cref{ass:abc_condition}. We therefore couple the flat regions with increasingly separated noise scales, distributing a finite variance budget among rare shocks. Every predictable stepsize policy then either makes insufficient progress in the flat coordinate or encounters damaging shocks whose conditional probabilities have a divergent sum, allowing us to apply the Borel--Cantelli lemma. The careful allocation of variance across scales is what reduces the loss relative to the critical exponents from an arbitrary polynomial factor to logarithmic factors. 

\begin{proof}[Theorem~\ref{thm:smooth_lb}~\ref{itm:smooth_lb_slow}]
Fix $L>0$. We start constructing an instance satisfying Assumptions~\ref{ass:convexity}-\ref{ass:abc_condition} which is independent of the stepsize policy $(\Gamma_k)_{k \geq 1}$ and $p$.
\paragraph{Objective.} For every integer $j\ge1$, define
\[
\varrho_j=2^{-3j},\qquad
\beta_j=\frac{\varrho_j}{4}\exp(-2^{j^2}),\qquad
R_j = 4 \exp(2^{j^2}),\qquad
T_j=\lceil R_j^4\rceil.
\]
The scale $R_j$ quantifies the cumulative stepsize needed to reach spatial scale $\varrho_j$, and $T_j$ is the iteration at which we test this progress. To construct the flat coordinate, define $\psi\colon[0,\infty)\to[0,\infty)$ by
\[
\psi(w)=
\begin{cases}
0, & w=0,\\
\beta_1, & w\ge 2\varrho_1,\\
\beta_j, & w\in[\varrho_j/2,2\varrho_j],\quad j\ge1,\\
\beta_{j+1}+(\beta_j-\beta_{j+1}) \chi\!\left((4/\varrho_j) \cdot w - 1\right), & w\in[\varrho_j/4,\varrho_j/2],\quad j\ge1,
\end{cases}
\]
where $\chi(t) = 3 \cdot t^2 - 2 \cdot t^3$. Now, set
\[
\varphi(v)=\int_0^{|v|}\psi(s)\,ds,\qquad v\in\mathbb R.
\]
Next, define 
\[
\chi_{\mathrm{b}}(u)=
\begin{cases}
-1/2, & u\le-1,\\
u+u^2/2, & -1<u<0,\\
u, & u\ge0.
\end{cases}
\]
Hence, let $\gamma_{\mathrm{b}}=\min\{1/4,L\}$ and define $F \colon \bbR^2 \rightarrow \bbR$ as 
\[
F(u,v) = \gamma_{\mathrm{b}} \cdot \bigl(\chi_{\mathrm{b}}(u)+\varphi(v)\bigr).
\]
Finally, set $\cX = [0,1]\times[-2,2]$ and $x_1=(0,2)$. By Lemma~\ref{lem:boundary-regular}, $F$ is convex and globally $L$-smooth. Its unique constrained minimizer is $(0,0)$, with $F_\star=0$ and $\nabla F(0,0)=(\gamma_{\mathrm{b}},0)\ne0$.
\paragraph{First-Order Oracle.} The role of the oracle is twofold: create rare but sufficiently large shocks in the $u$-coordinate, while keeping the second moment finite, and to add a benign unbiased perturbation in the $v$-coordinate. To this end, we use two independent noise variables.

\noindent First, let $Y_{\mathrm b}$ be a random variable with $Y_{\mathrm b}\ge1$ almost surely such that $\bbE[Y_{\mathrm b}] \le 2$ and 
\[
\forall j \ge 1,\ \forall 1 \le y \le T_j, \quad \bbP(Y_{\mathrm b}\ge y)\ge \vartheta_j^{\mathrm b} \cdot y^{-3/4},
\]
where $\vartheta_j^{\mathrm b}=2^{-(j+2)}T_j^{-1/4}$.
Thus $Y_{\mathrm b}$ retains a finite first moment, but still places heavy tails at each prescribed scale $T_j$. Such a random variable is constructed in the proof of Lemma~\ref{lem:boundary-mixture} by a suitable multiscale mixture of truncated heavy tails.

\noindent Next, define $\zeta \coloneqq (\Xi/2) \cdot \sqrt{Y_{\mathrm b}}$ where $\Xi$ is an independent Rademacher variable. Then $\zeta$ is centered, has finite variance, and satisfies
\[
\mathbb E[\zeta]=0, \qquad \mathbb E[\zeta^2]\le \frac12, \qquad |\zeta|\ge \frac12 \quad\text{a.s.}
\]
The lower-tail estimate on $Y_{\mathrm b}$ will later ensure that shocks of the order $\Gamma_k \cdot |\zeta_k|\gtrsim r_k^{-1}$ occur often enough.

\noindent Finally, let $\xi$ be independent of $(Y_{\mathrm b},\Xi)$ and uniformly distributed on $[1/2,3/2]$. At iteration $k$, we draw a copy $(Y_{\mathrm b,k},\Xi_k,\xi_k)$ of $(Y_{\mathrm b},\Xi,\xi)$ that is independent of $\cF_k$ and included in $\cF_{k+1}$, set $\zeta_k=(\Xi_k/2)\sqrt{Y_{\mathrm b,k}}$, and at $X_k=(u_k,v_k)$ the oracle returns
\[
\widehat\nabla_k F(X_k) = 
\begin{pmatrix}
\gamma_{\mathrm b}\cdot\chi_{\mathrm b}'(u_k)+\zeta_k\\[1mm]
\gamma_{\mathrm b}\cdot\xi_k\cdot\varphi'(v_k)
\end{pmatrix}.
\]
Since $\mathbb E[\xi]=1$ and $\mathbb E[\xi^2]=13/12$, we obtain
\[
\mathbb E_k[\widehat\nabla_k F(X_k)]=\nabla F(X_k), \qquad \mathbb E_k\bigl[\|\widehat\nabla_k F(X_k)\|^2\bigr] \le \frac12+\frac{25}{12} \cdot \gamma_{\mathrm b}^2<2.
\]
Hence Assumption~\ref{ass:abc_condition} holds with $a=b=0$ and $c=2$.

\paragraph{Error decomposition.}
Fix an arbitrary policy $(\Gamma_k)_{k\ge1}\in\cS_{\mathrm{fin}}$ and $p>0$. On the feasible set, $\chi_{\mathrm b}(u)=u$, so
\[
\Delta(X_k) = \gamma_{\mathrm b} \cdot \bigl(u_k+\varphi(v_k)\bigr),
\]
where both terms are nonnegative. The coordinate updates are
\begin{align*}
u_{k+1} &=\Pi_{[0,1]} \bigl(u_k-\Gamma_k \cdot (\gamma_{\mathrm b}+\zeta_k)\bigr),\\
v_{k+1} &=\Pi_{[-2,2]} \bigl(v_k-\gamma_{\mathrm b} \cdot \Gamma_k \cdot \xi_k \cdot \varphi'(v_k)\bigr).
\end{align*}
Set
\[
A_N=\sum_{k=1}^{N-1}\Gamma_k, \qquad r_k = (\log(e+k))^p \cdot k^{1/4},
\]
so that $1/r_k$ is the rate to be excluded. For integers $M,k\ge1$, define
\[
E_k^{(M)} \coloneqq \left\{\Gamma_k\zeta_k\le-\frac{2M}{r_k}\right\}, \qquad Q_M \coloneqq \sum_{k\ge1}\bbP_k\bigl(E_k^{(M)}\bigr).
\]
These events select negative shocks large enough, after multiplication by the stepsize, to force $u_{k+1}\ge M/r_k$ for sufficiently large $k$. Their conditional probabilities either have finite sum for some $M$,
or infinite sum for every $M$. We show that the first case leaves insufficient progress in the flat coordinate, whereas the second produces recurrent large errors in the noisy coordinate.

\paragraph{Control of the flat coordinate.} 
Suppose that $Q_M<\infty$ for some integer $M\ge1$. Lemma~\ref{lem:boundary-finite-hazard} gives
\[
\Gamma_k\longrightarrow0, \qquad A_{T_j}=o(R_j).
\]
Thus making the conditional probabilities of damaging shocks summable restricts the cumulative stepsize available for movement along $v$. The relevant obstruction is not merely that the stepsizes vanish, but that their sum is too small compared with $R_j$. To quantify the progress needed to approach zero, define
\[
H_2(w) \coloneqq \int_w^2\frac{ds}{\psi(s)}, \qquad 0<w\le2.
\]
This integral weights distance by the reciprocal slope: crossing a flatter region requires more cumulative stepsize. By Lemma~\ref{lem:sparse-profile}, $H_2$ is finite and strictly decreasing. Since $\psi=\beta_j$ on
$[\varrho_j,2\varrho_j]$,
\[
H_2(\varrho_j) \ge\int_{\varrho_j}^{2\varrho_j}\frac{ds}{\beta_j} =\frac{\varrho_j}{\beta_j} =R_j.
\]
Hence reaching $|v|\le\varrho_j$ requires $H_2(|v|)\ge R_j$. On the other hand, Lemma~\ref{lem:boundary-clock} ensures that $v_N\ne0$ for every $N$ and gives $H_2(|v_N|)\le C_\omega + 3 \cdot \gamma_{\mathrm b} \cdot A_N$ for all sufficiently large $N$, with a finite path-dependent constant $C_\omega\ge0$. Consequently, $H_2(|v_{T_j}|) \le C_\omega + 3 \cdot \gamma_{\mathrm b} \cdot A_{T_j} = o(R_j)$.
Comparing the required progress with the available progress and using the strict monotonicity of $H_2$, we obtain $|v_{T_j}|>\varrho_j$ eventually. The same constant-slope region also supplies a lower bound on the remaining objective value:
\[
\varphi(\varrho_j) \ge \int_{\varrho_j/2}^{\varrho_j}\beta_j\,ds = \frac{\varrho_j\cdot\beta_j}{2} = \frac{\varrho_j^2}{2R_j}.
\]
Since $\varphi$ is even and nondecreasing on $[0,\infty)$, and $T_j^{1/4}\ge R_j$, it follows that
\[
\begin{aligned}
r_{T_j} \cdot \Delta(X_{T_j}) &\ge \gamma_{\mathrm b} \cdot r_{T_j} \cdot \varphi(\varrho_j) \ge \frac{\gamma_{\mathrm b}}{2} \cdot \varrho_j^2 \cdot [\log(e+T_j)]^p\\
&\ge \frac{4^p \cdot \gamma_{\mathrm b}}{2} \cdot 2^{p\cdot j^2 - 6 \cdot j} \longrightarrow\infty, \qquad j\to\infty,
\end{aligned}
\]
where the last inequality follows from \eqref{eq:boundary-scale-growth}. The objective bound loses a factor $\varrho_j^2$ relative to $T_j^{-1/4}$; the chosen scales ensure that every positive power of $\log T_j$ dominates this loss.

\paragraph{Control of the noisy coordinate.}
Suppose now that $Q_M=\infty$ for every integer $M\ge1$. On $E_k^{(M)}$, the noise is negative and $|\zeta_k| \ge 1/2 \ge 2 \cdot \gamma_{\mathrm b}$. Thus the true gradient offsets at most half of the inward displacement caused by the noise. Since $u_k\ge0$,
\[
\begin{aligned}
u_{k+1} &=\Pi_{[0,1]} \bigl(u_k+\Gamma_k \cdot (|\zeta_k|-\gamma_{\mathrm b})\bigr) \ge\min\{1,\Gamma_k \cdot |\zeta_k|/2\} \ge\min\{1,M/r_k\}.
\end{aligned}
\]
For each fixed $M$, we eventually have $r_k\ge M$. Hence, at all sufficiently late occurrences of $E_k^{(M)}$,
\[
r_{k+1} \cdot \Delta(X_{k+1}) \ge \gamma_{\mathrm b} \cdot r_{k+1} \cdot u_{k+1} \ge\gamma_{\mathrm b} \cdot M \cdot \frac{r_{k+1}}{r_k} \ge\gamma_{\mathrm b} \cdot M,
\]
where we used that $r_k$ is increasing. Since $E_k^{(M)}\in\mathcal F_{k+1}$, conditional Borel--Cantelli (Lemma~\ref{lem:cbc}) gives infinitely many occurrences of $E_k^{(M)}$ on $\{Q_M=\infty\}$. Taking a countable intersection over $M$, the normalized error therefore exceeds arbitrarily large levels infinitely often, so that $\limsup_{k\to\infty} r_k \cdot \Delta(X_k)=\infty$.

\noindent Combining the two cases proves the failure of the claimed almost-sure convergence rate. The instance is independent of $p$. Intersecting these probability-one events for $p=1/m$, $m\ge1$, gives the conclusion simultaneously for all $p>0$: for any such $p$, choose $m$ with $1/m<p$ and use $[\log(e+k)]^p\ge[\log(e+k)]^{1/m}$. This probability-one event may depend on the fixed policy.
\end{proof}

\section{Nonsmooth Optimization}\label{section: main results of nonsmooth case}
In this section, we present our results on nonsmooth objectives. We start by specifying the first-order oracle, which complies with the following standard model.
\begin{assumption}\label{ass:oracle} 
There exists $\sigma>0$ such that, for every $k\ge 1$, it holds
 $$\bbE_k[\widehat{\nabla}_k F(X_k)] \in \partial F(X_k) \quad \text{and} \quad \bbE_k[\|\widehat{\nabla}_k F(X_k)\|^2]\le \sigma^2.$$
\end{assumption}
%

\noindent We first provide the lower bound on almost sure last-iterate convergence rates for SGD.
\begin{theorem}
\label{thm:nonsmooth_lb}
There exists a compact convex set $\mathcal X=[-1,1]\times[-2,2]\subset\bbR^2$, a nonsmooth function $F$ satisfying Assumptions \ref{ass:convexity}-\ref{ass:domain}, a first-order oracle satisfying Assumption \ref{ass:oracle} and $x_1\in\cX$ such that for every $(\Gamma_k)_{k \geq 1}\in\cS_{\mathrm{fin}}$ and every $p>0$, the iterates generated by eq.~\eqref{eq:sgd} satisfy $\Delta(X_k) \notin \fo_{a.s.}\left((\log(e+k))^{-p}\cdot k^{- 1/4}\right)$.
\end{theorem}
\noindent Unlike the $\widetilde{\mathcal O}(k^{-1/2})$ expected last-iterate rate \citep{Shamir2013}, an almost-sure guarantee under second moment assumptions must allow for slower polynomial decay in the worst case. The distinction between controlling the expected error and its decay along a sample path is therefore polynomial, not merely logarithmic.

\noindent We next provide the upper bound on almost sure last-iterate rates of SGD for anytime $k\ge 1$.
\begin{theorem}\label{thm:nonsmooth_ub}
Suppose that Assumptions \ref{ass:convexity}, \ref{ass:domain}, and \ref{ass:oracle} hold. Fix $\eta>0$. There exist $\alpha, t_0> 0$ and an $\alpha$-bounded predictable step-size sequence $(\Gamma_k)_{k \geq 1}$ such that the iterates generated by eq.~\eqref{eq:sgd} satisfy $\Delta(X_k) \in \fo_{a.s.}\bigl((\log(k+t_0))^{1/2+\eta} \cdot (k+t_0)^{-1/4}\bigr)$.
\end{theorem}
\noindent The proof uses deterministic, horizon-independent stepsizes. Together with \Cref{thm:nonsmooth_lb}, this shows that neither adaptivity nor randomization within $\cS_{\mathrm{fin}}$ improves the worst-case polynomial exponent over a suitable deterministic schedule. The optimal logarithmic correction remains open.


\subsection{Analysis}\label{sec:nonsmooth-sketch}
Here, we provide the proof of Theorem \ref{thm:nonsmooth_ub}. We defer to Appendix \ref{app:nonsmooth-endpoint} the proof of Theorem \ref{thm:nonsmooth_lb} that, instead, is similar in spirit to the argument described in \Cref{sec:smooth-sketch}.

\paragraph{Overview.} Existing almost sure convergence analysis of the last iterate \citep{Liu2024a} relies on smoothness and the absence of constraints, under which each step decreases the objective in expectation up to an error of order $\Gamma_k^2$. Without these conditions, a single step can increase the objective by an amount of order $\Gamma_k$, so the objective gap cannot be tracked directly by such a one-step supermartingale recursion. We instead track a smoothed surrogate of the constrained objective, in fact its \emph{Moreau envelope}, which nearly decreases in expectation from step to step. We note that this introduces a trade-off: stronger smoothing makes the surrogate deviate further from the objective, while weaker smoothing amplifies the per-step stochastic errors. Our key new ingredient resolves this trade-off by combining a time-varying envelope, whose smoothing vanishes at a carefully chosen speed, with a normalized supermartingale argument, yielding the optimal exponent $-1/4$.

\Needspace{5\baselineskip}
\begin{proof} Fix $\eta>0$, let $p=1/2+\eta$, fix deterministic $\gamma_0,\lambda_0>0$, choose $t_0\ge2$ with $\log(1+t_0)>4p$, and take the following stepsize schedule $\Gamma_k$ and smoothing parameter $\lambda_k$,
\[
 \Gamma_k=\gamma_0\cdot(k+t_0)^{-3/4}\cdot[\log(k+t_0)]^{-p},\ \
 \lambda_k=\lambda_0\cdot(k+t_0)^{-1/4}\cdot[\log(k+t_0)]^{-p},\ \
 r_k=\frac{(k+t_0)^{1/4}}{[\log(k+t_0)]^{p}}.
\]
Here $r_k$ is the reciprocal of the rate $[\log(k+t_0)]^p(k+t_0)^{-1/4}$. 
Note that $(\Gamma_k)_{k\ge1} \in \cS_\alpha$ for $\alpha = \Gamma_1$.

\paragraph{Moreau envelope.} 
Let $\Phi=\operatorname{cl}(F+\iota_\cX)$, the function whose epigraph is the closure of that of $F+\iota_\cX$. By Lemma~\ref{lem:closed-envelope}, $\Phi$ is proper, lower semicontinuous, convex, and such that $\min\Phi=F_\star$. Next, given a \emph{smoothing parameter} $\lambda>0$, the Moreau envelope of $\Phi$ \citep{BauschkeCombettes2017} is defined as
\[
\Phi_\lambda(x) \coloneqq \min_{y\in\cH}\Bigl\{\Phi(y)+\frac{\norm{y-x}^2}{2\lambda}\Bigr\},\qquad x\in\cH,
\]
and the minimization problem in \(y\) has a unique solution, denoted $p_\lambda(x)$, which is the proximal point of $x$ \citep{BauschkeCombettes2017}. By Lemma~\ref{lem:envelope}, $p_\lambda$ is well defined and $\Phi_\lambda$ is continuous.

\noindent Fix $x_\star\in\cX_\star$, which exists by Assumption~\ref{ass:domain}. Since $X_1=x_1\in\cX$ and every update ends with a projection onto $\cX$, all iterates lie in $\cX$. By Assumption~\ref{ass:oracle}, $\bbE_k[\widehat\nabla_kF(X_k)]\in\partial F(X_k)$; conditional Jensen gives $\bigl\|\bbE_k[\widehat\nabla_kF(X_k)]\bigr\|\le\sigma$, and the subgradient inequality gives
\begin{equation}\label{eq:nonsmooth-gap-integrability}
0\le\Delta(X_k)\le\ip{\bbE_k[\widehat\nabla_kF(X_k)]}{X_k-x_\star}\le\sigma\cdot\norm{X_k-x_\star}.
\end{equation}
Due to $x_\star \in \cX$, the nonexpansiveness of the projection, conditional unbiasedness, and Assumption~\ref{ass:oracle}, it holds that
\begin{equation}\label{eq:nonsmooth-distance-upper}
\bbE_k\bigl[\norm{X_{k+1}-x_\star}^2\bigr]\le \norm{X_k-x_\star}^2-2\cdot\Gamma_k\cdot\Delta(X_k)+\sigma^2\cdot\Gamma_k^2.
\end{equation}
Since $X_1=x_1$ is deterministic, $\bbE\bigl[\norm{X_1-x_\star}^2\bigr]=\norm{x_1-x_\star}^2<\infty$. Moreover, since $\Gamma_k\ge0$ and $\Delta(X_k)\ge0$, dropping the descent term in~\eqref{eq:nonsmooth-distance-upper} gives
\[
 \bbE_k\bigl[\norm{X_{k+1}-x_\star}^2\bigr]\le\norm{X_k-x_\star}^2+\sigma^2\cdot\Gamma_k^2 .
\]
Taking expectations, which is legitimate for nonnegative random variables by the tower property, yields $\bbE\bigl[\norm{X_{k+1}-x_\star}^2\bigr]\le\bbE\bigl[\norm{X_k-x_\star}^2\bigr]+\sigma^2\cdot\Gamma_k^2$. By induction,
\[
 \bbE\bigl[\norm{X_k-x_\star}^2\bigr]\le\norm{x_1-x_\star}^2+\sigma^2\cdot\sum_{i<k}\Gamma_i^2<\infty\qquad\text{for every }k,
\]
so $\Delta(X_k)$ is integrable by~\eqref{eq:nonsmooth-gap-integrability}. Taking full expectations in~\eqref{eq:nonsmooth-distance-upper}, summing over $k\le K$, dropping $\bbE\bigl[\norm{X_{K+1}-x_\star}^2\bigr]\ge0$, letting $K\to\infty$, and using $\sum_k\Gamma_k^2<\infty$ (Lemma~\ref{lem:quarter-schedule}) yield
\begin{equation}\label{eq:nonsmooth-energy-upper}
 2\cdot{\textstyle \sum_k}\Gamma_k\cdot\bbE[\Delta(X_k)]\le \norm{X_1-x_\star}^2+\sigma^2\cdot{\textstyle\sum_k}\Gamma_k^2<\infty.
\end{equation}

\noindent Set $Y_k=\Phi_{\lambda_k}(X_k)-F_\star$. By Lemma~\ref{lem:envelope}, $Y_k$ is $\cF_k$-measurable. By Lemma~\ref{lem:closed-envelope}, it holds that $\Phi(X_k)=F(X_k)$ and $\bbE_k[\widehat\nabla_kF(X_k)]\in\partial\Phi(X_k)$ for every $k$. Since $\Phi_{\lambda_k}\ge\min\Phi=F_\star$ and, by Lemma~\ref{lem:envelope-estimates}, $\Phi_{\lambda_k}\le F$ at feasible points, we have $0\le Y_k\le\Delta(X_k)$; thus, $Y_k$ is integrable. 

\noindent Since $\bbE_k[\|\widehat{\nabla}_kF(X_k)\|^2]\le\sigma^2$ and $\bigl\|\bbE_k[\widehat\nabla_kF(X_k)]\bigr\|\le\sigma$,
\[
 0\le\frac{\Gamma_k^2}{2\lambda_{k+1}}\cdot\bbE_k[\|{\widehat{\nabla}_k F(X_k)}\|^2]
        +\frac{\lambda_k-\lambda_{k+1}}2\cdot\bigl\|\bbE_k[\widehat\nabla_kF(X_k)]\bigr\|^2\le\frac{\sigma^2}{2}\cdot\Bigl(\frac{\Gamma_k^2}{\lambda_{k+1}}+\lambda_k-\lambda_{k+1}\Bigr).
\]
\noindent Since $\Delta(X_k)$ is integrable for every $k$ and $0\le\Phi_\lambda(X_k)-F_\star\le\Delta(X_k)$ for every $\lambda>0$, $\bbE_k[\|\widehat\nabla_kF(X_k)\|^2]\le\sigma^2$, $\Gamma_k$ is deterministic, and $\lambda_k$ is deterministic and nonincreasing (Lemma~\ref{lem:quarter-schedule}), Lemma~\ref{lem:envelope-step} gives

\begin{equation}\label{eq:recursion in main body}
   \bbE_k[Y_{k+1}]\le Y_k+\frac{\Gamma_k^2}{2\lambda_{k+1}}\cdot\bbE_k[\|{\widehat{\nabla}_k F(X_k)}\|^2]
        +\frac{\lambda_k-\lambda_{k+1}}2\cdot\bigl\|\bbE_k[\widehat\nabla_kF(X_k)]\bigr\|^2. 
\end{equation} 

\paragraph{Normalized supermartingale.} By Lemma~\ref{lem:quarter-schedule}, $r_k$ is nondecreasing with $r_{k+1}-r_k\le\Gamma_k/(4\cdot\gamma_0)$, and $\sum_k\Gamma_k/r_k=\sum_k\gamma_0/(k+t_0)=\infty$. Also note that Lemma~\ref{lem:quarter-schedule} gives 

$$\sum_kr_{k+1}\cdot\bbE\left[\frac{\Gamma_k^2}{2\lambda_{k+1}}\cdot\bbE_k[\|{\widehat{\nabla}_k F(X_k)}\|^2]
        +\frac{\lambda_k-\lambda_{k+1}}2\cdot\bigl\|\bbE_k[\widehat\nabla_kF(X_k)]\bigr\|^2\right]<\infty.$$ 

\noindent Besides, the proven fact $Y_k\le\Delta(X_k)$ and Equation~\eqref{eq:nonsmooth-energy-upper} give 
$$\sum_{k=1}^\infty\frac{\Gamma_k}{r_k}\cdot \left(r_k\cdot\bbE[Y_k]\right)=\sum_{k=1}^\infty\Gamma_k\cdot\bbE[Y_k]<\infty.$$

\noindent Applying Lemma~\ref{lem:scaled} to the recursion \eqref{eq:recursion in main body} with 
\begin{align*}
 a_k&=0, \qquad \theta_k=\Gamma_k, \qquad  G_k=Y_k \qquad \text{and}\\
 B_k&=\frac{\Gamma_k^2}{2\lambda_{k+1}}\cdot\bbE_k[\|{\widehat{\nabla}_k F(X_k)}\|^2]
        +\frac{\lambda_k-\lambda_{k+1}}2\cdot\bigl\|\bbE_k[\widehat\nabla_kF(X_k)]\bigr\|^2,
\end{align*}
it follows that $r_k\cdot Y_k\to0$ almost surely.

\noindent Applying Lemma \ref{lem:envelope-estimates} with $s=\bbE_k[\widehat\nabla_kF(X_k)]$ and using the fact $r_k\cdot\lambda_k=\lambda_0\cdot[\log(k+t_0)]^{-2p}\to 0$ proven in Lemma~\ref{lem:quarter-schedule}, it holds that
\[
 0\le r_k\cdot\Delta(X_k)-r_k\cdot Y_k\le\frac{r_k\cdot\lambda_k}{2}\cdot\bigl\|\bbE_k[\widehat\nabla_kF(X_k)]\bigr\|^2\le\frac{\sigma^2\cdot\lambda_0}{2\cdot[\log(k+t_0)]^{2p}}\longrightarrow0, \ \ \text{as} \ \ k\to\infty.
\]
Hence $r_k\cdot\Delta(X_k)\to0$ almost surely, that is, $\Delta(X_k)\in\fo_{a.s.}\bigl((\log(k+t_0))^{p}\cdot(k+t_0)^{-1/4}\bigr)$. 
\end{proof}

\section{Conclusions and Future Directions}
We characterized, up to polylogarithmic factors, the almost-sure anytime convergence rates of the last iterate of projected SGD on convex objectives. For smooth objectives under the conditional abc condition, the optimal exponent is $-1/3$ when a stationary minimizer exists and $-1/4$ otherwise. For nonsmooth objectives with bounded conditional second moments, it is again $-1/4$. 

\noindent Several directions remain open. One is to extend these results to the \emph{interpolation regime}, where the stochastic gradients vanish at the minimizer and faster rates may be possible. Another is to determine whether stepsize rules outside our class, such as gradient clipping and AdaGrad-Norm, whose effective stepsizes depend on the current oracle response, can improve the rates achieved here.

\newpage
\bibliographystyle{plainnat}
\bibliography{references}

\newpage
\appendix

\section{Details on Predictable Step-Sizes}\label{app:details_predictable_ss}
Here we explain the information convention of
Definition~\ref{def:predictable_ss}.
An internal random seed may be sampled independently of the oracle innovations
at the outset and included in $\cF_1$. At round $k$, both $X_k$ and $\Gamma_k$
are selected using $\cF_k$; the new sample is revealed only afterwards.
For a fixed i.i.d. oracle, this guarantees independence of the fresh
innovation from its predictable multiplier, even though the resulting
responses at adaptive query points are not independent.

\noindent\textit{Step sizes in $\cS_\alpha$.} Membership requires only that $\Gamma_k$ be computed from information available before the $k$-th oracle call and that $0<\Gamma_k\le\alpha$. The class therefore contains the following families.
\begin{itemize}[leftmargin=0.7cm]
\item \emph{Deterministic schedules} $\Gamma_k\equiv\gamma_k$ with $0<\gamma_k\le\alpha$, such as constant, polynomially decaying ($\gamma_0k^{-a}$ with $\gamma_0\le\alpha$), and step-decay schedules, as well as the logarithmically corrected schedules of Theorems~\ref{thm:smooth_ub} and~\ref{thm:nonsmooth_ub}. Horizon-dependent schedules concatenated over a deterministic sequence of horizons, as in the doubling trick, are included as well.
\item \emph{Adaptive rules based on past information} and capped by $\alpha$, such as the \emph{delayed} AdaGrad-Norm rule
$\Gamma_k=\min\{\alpha,\gamma_0/(1+\sum_{i<k}\|\widehat{\nabla}_iF(X_i)\|^2)^{1/2}\}$ with $\gamma_0>0$; see the scalar adaptive stepsizes of \citet{Li2019}.
\item \emph{Randomized policies} whose internal randomization is independent of the current oracle sample and revealed before the oracle call, for instance $\Gamma_k=U_k\gamma_k$ with $U_k\in(0,1]$ drawn from the algorithmic seed, as in {random} learning rates \citep{Musso2020}.
\end{itemize}

\noindent\textit{Step sizes outside $\cS_\alpha$.} A rule leaves $\cS_\alpha$, and also $\cS_{\mathrm{fin}}$, as soon as its dependence on the fresh oracle response makes its effective stepsize non-\(\mathcal F_k\)-measurable. Written in the form~\eqref{eq:sgd} with the original oracle response, this is the case for the following rules.
\begin{itemize}[leftmargin=0.7cm]
\item \emph{Gradient clipping} \citep{Zhang2020,Gorbunov2020}, $\Gamma_k=\gamma_k\min\{1,\hat{c}_k/\|\widehat{\nabla}_kF(X_k)\|\}$, where $\hat{c}_k$ is the clipping level.
\item \emph{Normalized SGD} \citep{Hazan2015}, $\Gamma_k=\gamma_k/\|\widehat{\nabla}_kF(X_k)\|$.
\item \emph{The usual AdaGrad-Norm update} \citep{Ward2020}, $\Gamma_k=\gamma_0/(b_0^2+\sum_{i\le k}\|\widehat{\nabla}_iF(X_i)\|^2)^{1/2}$, whose delayed version above is predictable.
\item \emph{Stochastic Polyak stepsize} \citep{Loizou2021}, $\Gamma_k=\min\{\bar\gamma,(f_k(X_k)-f_k^\star)/(c_{\mathrm{P}}\|\widehat{\nabla}_kF(X_k)\|^2)\}$ with $c_{\mathrm{P}},\bar\gamma>0$, where $f_k$ is the current sample loss, $\nabla f_k(X_k)=\widehat{\nabla}_kF(X_k)$, and $f_k^\star=\inf f_k$; it extends the deterministic Polyak stepsize \citep{Polyak1987}.
\end{itemize}
Such rules can shrink exactly the steps whose oracle response is large, so the conditional probability of a damaging shock is no longer the survival function of the noise evaluated at a predictable threshold, which is the mechanism behind our lower bounds. Whether they can improve the almost-sure exponents is an open question. Thus timing, rather than determinism, separates the two classes. Finally, per-coordinate or matrix preconditioners, such as diagonal AdaGrad or Adam, are not scalar multipliers and lie outside the form~\eqref{eq:sgd} altogether.

\noindent The class $\cS_\alpha$ retains the strictly positive, deterministic-cap
convention of the main smooth theorems. The boundary and nonsmooth lower-bound
proofs hold for the larger class $\cS_{\mathrm{fin}}$, which has no deterministic cap and allows zero steps. The stationary hard curvature
uses the prescribed $\alpha$ and is never selected from a realized random
bound. 

\newpage
\section{Known Results and Numerical Estimates}\label{app:tools}

This appendix collects the known results and numerical estimates used in our analysis. Constructions and estimates specific to a single lower-bound instance are deferred to Appendices~\ref{app:smooth-lb} and~\ref{app:nonsmooth-endpoint}.

\noindent Throughout this appendix, $X_k$ and $\Gamma_k$ are
$\cF_k$-measurable, and a round-$k$ event belongs to $\cF_{k+1}$.
All inequalities between random variables below hold almost surely.



\label{app:common tools}

\paragraph{Probabilistic results.}
We begin with the standard conditional Borel-Cantelli Lemma and Robbins-Siegmund Theorem. Then we present a corollary of Robbins-Siegmund Theorem in Lemma \ref{lem:scaled}. 


\begin{lemma}[Conditional Borel--Cantelli, \cite{durrett2019probability}, Theorem~4.3.4]\label{lem:cbc}
Let $E_k\in\cF_{k+1}$ and $\pi_k=\bbP_k(E_k)$. Then
\begin{equation}\label{eq:cbc}
 \left\{\sum_{k\ge1}\pi_k=\infty\right\}
 \subseteq\{E_k\text{ occurs infinitely often}\}
 \quad\text{up to a null set}.
\end{equation}
\end{lemma}

\begin{theorem}[Robbins-Siegmund \citep{Robbins1971}]\label{lem:rs}
Let $G_k,D_k,B_k$ be nonnegative adapted integrable random variables,
and let $a_k\ge0$ be deterministic with $\sum_k a_k<\infty$. Suppose
\begin{equation}\label{eq:rs-recursion}
 \bbE_k[G_{k+1}]\le(1+a_k)G_k-D_k+B_k,
 \qquad \sum_k\bbE[B_k]<\infty.
\end{equation}
Then $G_k$ has a finite almost-sure limit and $\sum_kD_k<\infty$ almost
surely.
\end{theorem}

\begin{lemma}\label{lem:scaled}
Let $G_k\ge0$ be adapted and integrable, let $a_k\ge0$ be deterministic
and summable, and let $B_k\ge0$ be adapted. Suppose
\[
 \bbE_k[G_{k+1}]\le(1+a_k)G_k+B_k.
\]
Let $r_k,\theta_k>0$ be deterministic. Assume, after possibly discarding
a finite prefix, there exists a finite constant $C_r$ such that
\begin{equation}\label{eq:scaled-hypotheses}
 \begin{gathered}
 r_k\text{ is nondecreasing},\qquad
 r_{k+1}-r_k\le C_r\theta_k,\qquad
 \sum_k\theta_k\bbE[G_k]<\infty,\\
 \sum_k r_{k+1}\bbE[B_k]<\infty,\qquad
 \sum_k\frac{\theta_k}{r_k}=\infty.
 \end{gathered}
\end{equation}
Then $r_kG_k\to0$ almost surely.
\end{lemma}
\begin{proof}
Discard a deterministic finite prefix so that all eventual hypotheses hold;
this does not affect the conclusion. Set $Z_k=r_kG_k$. Its conditional recursion is
\[
 \bbE_k[Z_{k+1}]\le(1+a_k)Z_k+\beta_k,
 \quad
 \beta_k=(1+a_k)(r_{k+1}-r_k)G_k+r_{k+1}B_k.
\]
Because a summable nonnegative sequence is bounded,
$\sum_k\bbE[\beta_k]<\infty$ follows from~\eqref{eq:scaled-hypotheses}.
Theorem~\ref{lem:rs}, with zero descent term, gives a finite nonnegative
limit $Z_\infty$. Tonelli's theorem also gives
$\sum_k\theta_kG_k<\infty$ almost surely. On any path where
$Z_\infty>0$, eventually $Z_k\ge Z_\infty/2$, so
\[
 \sum_k\theta_kG_k
 =\sum_k\frac{\theta_k}{r_k}Z_k=\infty,
\]
a contradiction. Hence $Z_\infty=0$ almost surely.
\end{proof}

\paragraph{Convex Analysis.} We first collect the convex-analysis facts used below, starting from the following basic fact.
\begin{lemma}[\cite{BauschkeCombettes2017}, Chapters 9 and 16]\label{lem:closed-envelope}
Under Assumptions~\ref{ass:convexity} and~\ref{ass:domain}, define
\begin{equation}\label{eq:closed-constrained}
 \Phi:=\operatorname{cl}(F+\iota_\cX),
\end{equation}
where $\operatorname{cl}$ denotes the lower-semicontinuous convex closure
(epigraph closure). Then $\Phi$ is proper, lower semicontinuous, convex,
has domain contained in $\cX$, and $\min\Phi=F_\star$, attained at each
$x_\star\in\cX_\star$. Whenever $x\in\cX$ has a subgradient
$s\in\partial F(x)$,
\begin{equation}\label{eq:closure-values}
 \Phi(x)=F(x),\qquad s\in\partial\Phi(x).
\end{equation}
If $F$ is continuous, then $\Phi=F+\iota_\cX$.
\end{lemma}

\noindent Lemma~\ref{lem:envelope} collects standard properties of proximal points and \emph{Moreau envelopes} \citep{Moreau1965,BauschkeCombettes2017}.

\begin{lemma}[\cite{Moreau1965}; \cite{BauschkeCombettes2017}, Section~12.4]\label{lem:envelope}
Under Assumptions~\ref{ass:convexity} and~\ref{ass:domain}, let $\Phi$ be the closed constrained objective in~\eqref{eq:closed-constrained}. For $\lambda>0$ and $x\in\cH$, set
\begin{equation}\label{eq:envelope-def}
 \begin{aligned}
 \Phi_\lambda(x)&:=\min_{y\in\cH}
       \left\{\Phi(y)+\frac{\norm{y-x}^2}{2\lambda}\right\},\\
 p_\lambda(x)&:=\prox_{\lambda\Phi}(x),\qquad
 v_\lambda(x):=\frac{x-p_\lambda(x)}{\lambda}.
 \end{aligned}
\end{equation}
The minimizer exists uniquely, belongs to $\cX$, and
$v_\lambda(x)\in\partial\Phi(p_\lambda(x))$. The maps $p_\lambda$ and
$\Phi_\lambda$ are continuous.
\end{lemma}


\noindent Lemma~\ref{lem:envelope-estimates} adds the subgradient and parameter-variation estimates needed in our analysis, analogous to results in \citet{Brezis1973,BotBohm2020,BohmWright2021}. We include its proof for completeness.

\begin{lemma}\label{lem:envelope-estimates}
Let $\Phi_\lambda$, $p_\lambda$, and $v_\lambda$ be as in Lemma~\ref{lem:envelope}. If $x\in\cX$ and $s\in\partial F(x)$, then, for every $\lambda>0$,
\begin{align}
 \ip{s}{v_\lambda(x)}&\ge\norm{v_\lambda(x)}^2,
 &\norm{v_\lambda(x)}&\le\norm{s},\label{eq:envelope-monotone}\\
 0\le F(x)-\Phi_\lambda(x)&\le\frac\lambda2\norm{s}^2,
 &&\label{eq:envelope-gap}
\end{align}
and, for $0<\lambda'\le\lambda$,
\begin{equation}\label{eq:envelope-parameter}
 0\le\Phi_{\lambda'}(x)-\Phi_\lambda(x)
 \le\frac{\lambda-\lambda'}2\norm{s}^2.
\end{equation}
\end{lemma}
\begin{proof}
Lemma~\ref{lem:closed-envelope} puts $s$ in $\partial\Phi(x)$ and
identifies $\Phi(x)=F(x)$. Since $v_\lambda(x)\in\partial\Phi(p_\lambda(x))$ by Lemma~\ref{lem:envelope}, monotonicity of $\partial\Phi$ between $x$ and $p_\lambda(x)$
yields
$\ip{s-v_\lambda(x)}{x-p_\lambda(x)}\ge0$. Dividing by $\lambda$ gives the first inequality in~\eqref{eq:envelope-monotone}, and the Cauchy--Schwarz inequality then gives the second. The competitor $y=x$ in~\eqref{eq:envelope-def} proves
$\Phi_\lambda(x)\le\Phi(x)=F(x)$. The affine support at $x$ also gives
\[
 \Phi_\lambda(x)\ge F(x)+\inf_y
 \left\{\ip{s}{y-x}+\frac{\norm{y-x}^2}{2\lambda}\right\}
 =F(x)-\frac\lambda2\norm{s}^2,
\]
proving~\eqref{eq:envelope-gap}.
For~\eqref{eq:envelope-parameter}, the lower inequality follows from~\eqref{eq:envelope-def}, since the penalty $\norm{y-x}^2/(2\lambda)$ is nonincreasing in $\lambda$. Put $z=p_{\lambda'}(x)$ and
$v=v_{\lambda'}(x)\in\partial\Phi(z)$. The support at $z$ gives
\[
 \Phi_\lambda(x)\ge\Phi(z)+\ip{v}{x-z}-\frac\lambda2\norm v^2
 =\Phi_{\lambda'}(x)-\frac{\lambda-\lambda'}2\norm v^2,
\]
where the equality uses $x-z=\lambda'v$ and $\Phi_{\lambda'}(x)=\Phi(z)+\lambda'\norm v^2/2$.
Since $\norm v\le\norm s$ by~\eqref{eq:envelope-monotone} with $\lambda'$ in place of $\lambda$, this proves the upper inequality.
\end{proof}

\noindent Below, Equation~\eqref{eq:fixed-envelope-step} is the convex case of the one-step estimate of \citet{DavisDrusvyatskiy2019}, and~\eqref{eq:variable-envelope-step} extends it to decreasing smoothing parameters \citep{BotBohm2020,BohmWright2021}. We include the proofs here for completeness.
\begin{lemma}\label{lem:envelope-step}
Let $(X_k)_{k\ge1}$ be generated by~\eqref{eq:sgd}, let $\Phi_\lambda$ and $v_\lambda$ be as in Lemma~\ref{lem:envelope}, and suppose that $\bbE_k[\widehat{\nabla}_k F(X_k)]\in\partial F(X_k)$ and that $\Phi_\lambda(X_k)$, $\Phi_\lambda(X_{k+1})$, $\|\widehat{\nabla}_kF(X_k)\|^2$, and $\Gamma_k^2\|\widehat{\nabla}_kF(X_k)\|^2$ are integrable for any $\lambda>0$. For $\lambda>0$ and a predictable
$\Gamma_k\ge0$,
\begin{equation}\label{eq:fixed-envelope-step}
 \bbE_k[\Phi_\lambda(X_{k+1})-F_\star]
 \le\Phi_\lambda(X_k)-F_\star
      -\Gamma_k\norm{v_\lambda(X_k)}^2
      +\frac{\Gamma_k^2}{2\lambda}\bbE_k[\|{\widehat{\nabla}_k F(X_k)}\|^2].
\end{equation}
For positive deterministic nonincreasing $\lambda_k$, set
$Y_k=\Phi_{\lambda_k}(X_k)-F_\star$. Then $0\le Y_k\le\Delta(X_k)$ and
\begin{equation}\label{eq:variable-envelope-step}
 \bbE_k[Y_{k+1}]\le Y_k+\frac{\Gamma_k^2}{2\lambda_{k+1}}\bbE_k[\|{\widehat{\nabla}_k F(X_k)}\|^2]
        +\frac{\lambda_k-\lambda_{k+1}}2\norm{\bbE_k[\widehat{\nabla}_k F(X_k)]}^2.
\end{equation}
\end{lemma}
\begin{proof}
Put $z=p_\lambda(X_k)\in\cX$ and $v=(X_k-z)/\lambda$. The competitor
$z$ in the envelope at $X_{k+1}$ and projection nonexpansiveness imply
\[
 \begin{split}
 \Phi_\lambda(X_{k+1})
 &\le\Phi(z)+\frac{\norm{X_{k+1}-z}^2}{2\lambda}\\
 &\le\Phi(z)+\frac{\norm{X_k-\Gamma_k\widehat{\nabla}_k F(X_k)-z}^2}{2\lambda}\\
 &=\Phi_\lambda(X_k)-\Gamma_k\ip{\widehat{\nabla}_k F(X_k)}{v}
                    +\frac{\Gamma_k^2}{2\lambda}\norm{\widehat{\nabla}_k F(X_k)}^2.
 \end{split}
\]
Take conditional expectations and apply~\eqref{eq:envelope-monotone} of Lemma~\ref{lem:envelope-estimates}
with $s=\bbE_k[\widehat{\nabla}_k F(X_k)]$ to obtain~\eqref{eq:fixed-envelope-step}.
Apply it with $\lambda=\lambda_{k+1}$, drop its nonpositive term, and
use~\eqref{eq:envelope-parameter} at $X_k$ to prove
\eqref{eq:variable-envelope-step}. The lower bound on $Y_k$ follows from
$\min\Phi=F_\star$ and the upper bound from the competitor $X_k$.
\end{proof}

\paragraph{Numerical Estimates.} We close with four important estimates that will be used in later proofs.
\begin{lemma}\label{lem:power-sums}
If $0\le a_0<1$ and $b_0\ge0$, then a finite deterministic constant
$C_{a_0,b_0}$ exists such that
\begin{equation}\label{eq:weighted-power-sum}
 \sum_{k=1}^{N-1}k^{-a_0}[\log(e+k)]^{-b_0}
 \le C_{a_0,b_0}N^{1-a_0}[\log(e+N)]^{-b_0},\qquad N\ge2.
\end{equation}
\end{lemma}
\begin{proof}
For $k\le\sqrt N$, ignore the logarithmic factor, which is at most one. The resulting sum is at most $C N^{(1-a_0)/2}$. For
$\sqrt N<k<N$, the inequality $\log(e+k)\ge c\log(e+N)$ holds with a universal $c>0$ for all large $N$. Summing $k^{-a_0}$ on that range gives at most $C'N^{1-a_0}[\log(e+N)]^{-b_0}$. The first term is smaller than the latter bound eventually, because
$N^{-(1-a_0)/2}[\log(e+N)]^{b_0}\to0$. Increasing the constant handles finitely many smaller $N$.
\end{proof}


\begin{lemma}\label{lem:clock}
Let $h:[0,\bar b]\to[0,\infty)$ be nondecreasing and $K_0$-Lipschitz, with
$h(0)=0$ and $h(w)>0$ for $w>0$. Let
$\mathscr C(w)=\int_w^{\bar b} ds/h(s)$ for $0<w\le \bar b$. If
$0<w\le \bar b$, $\epsilon\ge0$, and $\epsilon K_0\le1/2$, put
$w^+=w-\epsilon h(w)$. Then
\begin{equation}\label{eq:clock-step}
 w/2\le w^+\le w,\qquad
 \mathscr C(w^+)-\mathscr C(w)\le2\epsilon.
\end{equation}
\end{lemma}
\begin{proof}
Since $h(w)\le K_0w$, the bound on $w^+$ follows, and in particular
$w^+>0$. Moreover,
\[
 h(w^+)\ge h(w)-K_0(w-w^+)
             =(1-\epsilon K_0)h(w)\ge h(w)/2.
\]
Monotonicity gives
$\int_{w^+}^w ds/h(s)\le(w-w^+)/h(w^+)\le2\epsilon$.
The case $\epsilon=0$ is included. No continuous-time approximation is
used.
\end{proof}

\begin{lemma}\label{lem:sparse-profile}
Let $\varrho_j$, $\beta_j$, $R_j$, $T_j$, and $\psi$ be as in Section~\ref{sec:smooth-sketch}, so that $R_j=\varrho_j/\beta_j$. The function $\psi$ belongs to $C^1([0,\infty))$, is nondecreasing, and
satisfies
\begin{equation}\label{eq:profile-bounds}
 \psi'(0)=0,\qquad 0\le\psi'(w)\le1,\qquad
 0<\psi(w)\le\min\{w,\beta_1\}\quad(w>0).
\end{equation}
The even potential
\begin{equation}\label{eq:shared-potential}
 \varphi(v):=\int_0^{|v|}\psi(s)\,ds
\end{equation}
is convex and $C^2$, with $\varphi''(0)=0$,
$0\le\varphi''(v)\le1$, $|\varphi'(v)|\le\beta_1<1$, and
$\varphi(v)=0$ exactly when $v=0$.
For ${\bar{b}}\ge2\varrho_1=1/4$, the clock
\begin{equation}\label{eq:shared-clock}
 H_{\bar{b}}(w)=\int_w^{\bar{b}}\frac{ds}{\psi(s)},\qquad 0<w\le {\bar{b}},
\end{equation}
is finite and strictly decreasing, and
\begin{equation}\label{eq:plateau-estimates}
 \varphi(\varrho_j)\ge\frac{\varrho_j\beta_j}{2},\qquad
 H_{\bar{b}}(\varrho_j)\ge R_j.
\end{equation}
Finally, for every $p>0$,
\begin{equation}\label{eq:boundary-scale-growth}
 \begin{gathered}
 T_j^{1/4}\ge R_j,\qquad
 \frac{T_j^{1/4}}{R_j}\to1,\qquad
 \log(e+T_j)\ge4\cdot2^{j^2},\\
 \varrho_j^2\,[\log(e+T_j)]^{p}\ge4^p\,2^{pj^2-6j}\longrightarrow\infty, \ \ \text{as} \ \ j\to\infty.
 \end{gathered}
\end{equation}
\end{lemma}
\begin{proof}
The heights $\beta_j$ decrease with $j$. Since
$\chi'(t)=6t(1-t)$ is nonnegative, vanishes at both endpoints, and is at
most $3/2$, the profile is nondecreasing and $C^1$ away from zero. On a
transition interval $[\varrho_j/4,\varrho_j/2]$,
\[
 0\le\psi'(w)
 \le\frac{6(\beta_j-\beta_{j+1})}{\varrho_j}
 \le\frac32e^{-2^{j^2}}<1.
\]
On a plateau $[\varrho_j/2,2\varrho_j]$ its derivative is zero. The intervals
$[\varrho_j/4,2\varrho_j]$ cover $(0,2\varrho_1]$ and meet at their
endpoints. On each such interval,
\[
 0\le\frac{\psi(w)}w\le\frac{4\beta_j}{\varrho_j}
                  =e^{-2^{j^2}}.
\]
Thus $\psi(w)/w\to0$ at zero, so $\psi$ is differentiable there with
derivative zero. The displayed derivative bounds also tend to zero as
$w\downarrow0$, proving continuity of $\psi'$ at zero. Its global
$1$-Lipschitz property and $\psi(0)=0$ give $\psi(w)\le w$; the maximum
height is $\beta_1$. Strict positivity follows from the positive heights
and their interpolations.

\noindent For $v\ne0$,
$\varphi'(v)=\sgn(v)\psi(|v|)$ and
$\varphi''(v)=\psi'(|v|)$. At zero,
$\varphi'(0)=0$ and
$\varphi''(0)=\lim_{w\downarrow0}\psi(w)/w=0$.
This proves $C^2$ regularity and the derivative bounds. Nonnegativity of
the second derivative proves convexity. The integral is strictly
positive off zero.

\noindent The clock is finite because $\psi$ is continuous and bounded away from
zero on each compact subinterval of $(0,\bar b]$. It is strictly decreasing
since its integrand is positive. Integrating over a plateau gives
\[
 \varphi(\varrho_j)\ge\int_{\varrho_j/2}^{\varrho_j}\beta_j\,ds,
 \qquad
 H_{\bar b}(\varrho_j)\ge\int_{\varrho_j}^{2\varrho_j}\frac{ds}{\beta_j},
\]
which proves~\eqref{eq:plateau-estimates}.

\noindent Finally, $R_j^4\le T_j<R_j^4+1$ gives the first two relations in~\eqref{eq:boundary-scale-growth}. Moreover, $\log(e+T_j)\ge4\log R_j=4(\log4+2^{j^2})\ge4\cdot2^{j^2}$, and $\varrho_j^2=2^{-6j}$ gives the last one.
\end{proof}

\begin{lemma}\label{lem:quarter-schedule}
Let $p>1/2$ and $\gamma_0,\lambda_0>0$, and let $t_0\ge2$ satisfy $\log(1+t_0)>4p$. Put
\begin{equation}\label{eq:quarter-scales}
 \begin{gathered}
 \Gamma_k=\gamma_0 (k+t_0)^{-3/4}[\log(k+t_0)]^{-p},\qquad
 \lambda_k=\lambda_0(k+t_0)^{-1/4}[\log(k+t_0)]^{-p},\\
 r_k=(k+t_0)^{1/4}[\log(k+t_0)]^{-p}.
 \end{gathered}
\end{equation}
Then $\sum\Gamma_k^2<\infty$, $\lambda_k$ decreases, $r_k$ increases,
and
\begin{align}
 &0\le r_{k+1}-r_k\le\frac{\Gamma_k}{4\gamma_0},
 &\frac{\Gamma_k}{r_k}&=\frac{\gamma_0}{k+t_0},\label{eq:quarter-increments}\\
 &\frac{\Gamma_k^2}{\lambda_{k+1}}+(\lambda_k-\lambda_{k+1})
 \le \frac{C}{(k+t_0)^{5/4}[\log(k+t_0)]^{p}},
 &r_k\lambda_k&=\frac{\lambda_0}{[\log(k+t_0)]^{2p}}\to0,\label{eq:quarter-errors}\\
 &\sum_k r_{k+1}
 \left[\frac{\Gamma_k^2}{\lambda_{k+1}}+
                   \lambda_k-\lambda_{k+1}\right]<\infty,
 &\sum_k\frac{\Gamma_k}{r_k}&=\infty.\label{eq:quarter-sums}
\end{align}
Here $C$ is finite and deterministic.
\end{lemma}
\begin{proof}
Since $\log(k+t_0)\ge\log(1+t_0)>0$, the stepsize squares are bounded by a
constant multiple of $(k+t_0)^{-3/2}$.
For $R(x)=x^{1/4}(\log x)^{-p}$,
\[
 R'(x)=x^{-3/4}(\log x)^{-p}
                   \left(\frac14-\frac p{\log x}\right).
\]
The bracket is positive for $x\ge1+t_0$, because $\log(1+t_0)>4p$, and is at most $1/4$.
The remaining factor decreases, so integration over $[k+t_0,k+1+t_0]$
proves the increment inequality. For
$J(x)=\lambda_0x^{-1/4}(\log x)^{-p}$,
\[
 -J'(x)=\lambda_0x^{-5/4}(\log x)^{-p}
                    \left(\frac14+\frac p{\log x}\right),
\]
which gives the claimed bound on $\lambda_k-\lambda_{k+1}$, since $p/\log x<1/4$ there.
Since $(k+1+t_0)/(k+t_0)$ and $\log(k+1+t_0)/\log(k+t_0)$ are uniformly bounded,
$\Gamma_k^2/\lambda_{k+1}\le C (k+t_0)^{-5/4}[\log(k+t_0)]^{-p}$ and
$r_{k+1}\le C r_k$. Hence the summand in \eqref{eq:quarter-sums} is at
most $C/\{(k+t_0)[\log(k+t_0)]^{2p}\}$, whose sum is finite by the integral test
because $2p>1$. The identities for $r_k\lambda_k$ and $\Gamma_k/r_k$
are immediate, and the latter has a divergent harmonic sum.
\end{proof}

\newpage
\section{Missing Proof of Theorem~\ref{thm:smooth_lb}}\label{app:smooth-lb}
This appendix completes the proof of Theorem~\ref{thm:smooth_lb}: Appendix~\ref{app:boundary-lb} provides the missing proof of Theorem \ref{thm:smooth_lb}~\ref{itm:smooth_lb_slow} in Section~\ref{sec:smooth-sketch}; Appendix~\ref{app:smooth-endpoint} proves Theorem \ref{thm:smooth_lb}~\ref{itm:smooth_lb_fast}.

\subsection{Missing Proof of Theorem~\ref{thm:smooth_lb}~\ref{itm:smooth_lb_slow}}
\label{app:boundary-lb}

This subsection proves the lemmas invoked in the proof of Theorem~\ref{thm:smooth_lb}~\ref{itm:smooth_lb_slow} in Section~\ref{sec:smooth-sketch}, where the objective, the oracle, and the quantities $A_N$, $r_k$, $E_k^{(M)}$, and $Q_M$ are defined. 


\begin{lemma}\label{lem:boundary-regular}
The objective $F(u,v)=\gamma_{\mathrm{b}}\bigl(\chi_{\mathrm{b}}(u)+\varphi(v)\bigr)$ of Section~\ref{sec:smooth-sketch} is convex, continuously
differentiable, bounded below, globally $1$-Lipschitz, and has globally
$L$-Lipschitz gradient. On $\cX=[0,1]\times[-2,2]$,
\begin{equation}\label{eq:boundary-linear-gap}
 F(u,v)-F_\star=\gamma_{\mathrm{b}}\bigl(u+\varphi(v)\bigr).
\end{equation}
Its unique constrained minimizer is $(0,0)$, with $F_\star=0$ and
$\nabla F(0,0)=(\gamma_{\mathrm{b}},0)\ne0$.
\end{lemma}
\begin{proof}
The derivative $\chi_{\mathrm{b}}'$ is the continuous function equal to
$0$, $1+u$, and $1$ on the three regions in the definition of
$\chi_{\mathrm{b}}$. It is nondecreasing, takes values in
$[0,1]$, and is $1$-Lipschitz. Thus $\chi_{\mathrm{b}}$ is convex,
$C^1$, and bounded below by $-1/2$. By
Lemma~\ref{lem:sparse-profile}, $\varphi$ is convex, $C^2$, nonnegative,
and has a $1$-Lipschitz derivative of magnitude at most $\beta_1<1$.
The separable gradient consequently has Lipschitz constant at most
$\gamma_{\mathrm{b}}\le L$, while
\[
 \norm{\nabla F(u,v)}
 \le\gamma_{\mathrm{b}}\sqrt{1+\beta_1^2}<1.
\]
The latter bound proves global $1$-Lipschitz continuity of $F$ by integration
along line segments. Also $F\ge-\gamma_{\mathrm{b}}/2$ on $\bbR^2$.
On the feasible set $\chi_{\mathrm{b}}(u)=u\ge0$; both terms vanish
simultaneously exactly at $(0,0)$. This proves the minimum claim and
\eqref{eq:boundary-linear-gap}. The derivative at that point follows
from $\chi_{\mathrm{b}}'(0)=1$ and $\varphi'(0)=0$.
\end{proof}

\begin{lemma}\label{lem:boundary-mixture}
Let $(T_j)_{j\ge1}$ be the sequence defined in Section~\ref{sec:smooth-sketch} and put $\delta_j=2^{-j-2}$. There exists a random variable $Y_{\mathrm b}$ satisfying $1\le Y_{\mathrm b}<\infty$ almost surely and
\begin{equation}\label{eq:boundary-mixture-estimates}
\bbE[Y_{\mathrm b}]\le2,\qquad \bbP(Y_{\mathrm b}\ge y)\ge\delta_jT_j^{-1/4}y^{-3/4}\qquad(j\ge1,\ 1\le y\le T_j).
\end{equation}
If $\Xi$ is independent of $Y_{\mathrm b}$ with $\bbP(\Xi=1)=\bbP(\Xi=-1)=1/2$, then $\zeta=\Xi\sqrt{Y_{\mathrm b}}/2$ satisfies
\[
\bbE[\zeta]=0,\qquad \bbE[\zeta^2]\le\frac12,\qquad |\zeta|\ge\frac12 \quad\text{almost surely}.
\]
Moreover, for every $j\ge1$,
\begin{equation}\label{eq:boundary-noise-tail}
\bbP(\zeta\le-t) \ge 2^{-j-3}T_j^{-1/4}(4t^2)^{-3/4}, \qquad \frac12\le t\le\frac{\sqrt{T_j}}2.
\end{equation}
\end{lemma}
\begin{proof}
Select a random variable $V_{\mathrm b}\ge1$ with
\[
\bbP(V_{\mathrm b}\ge y)=y^{-3/4}, \qquad y\ge1,
\]
and define
\[
\vartheta_j^{\mathrm b} = 2^{-j-2}T_j^{-1/4}, \qquad j\ge1.
\]
Since $T_j\ge1$,
\[
\sum_{j\ge1}\vartheta_j^{\mathrm b} \le\sum_{j\ge1}2^{-j-2} =\frac14.
\]
Thus we may choose an independent random variable $J\in\{0,1,2,\ldots\}$ with
\[
\bbP(J=j)=\vartheta_j^{\mathrm b}\quad(j\ge1), \qquad \bbP(J=0)=1-\sum_{j\ge1}\vartheta_j^{\mathrm b}.
\]
Set
\[
Y_{\mathrm b} = 
\begin{cases}
1, & J=0,\\
\min\{V_{\mathrm b},T_J\}, & J\ge1.
\end{cases}
\]
Then $1\le Y_{\mathrm b}<\infty$ almost surely.

\noindent The layer-cake formula gives
\[
\bbE[\min\{V_{\mathrm b},T_j\}] =1+\int_1^{T_j}y^{-3/4}\,dy =4T_j^{1/4}-3.
\]
Consequently, independence and nonnegative summation yield
\[
\begin{aligned}
\bbE[Y_{\mathrm b}] &=1-\sum_{j\ge1}\vartheta_j^{\mathrm b}+\sum_{j\ge1}\vartheta_j^{\mathrm b} (4T_j^{1/4}-3)\\
&=1+4\sum_{j\ge1}2^{-j-2} -4\sum_{j\ge1}\vartheta_j^{\mathrm b} \le2.
\end{aligned}
\]
For every $j\ge1$ and $1\le y\le T_j$,
\[
\begin{aligned}
\bbP(Y_{\mathrm b}\ge y) &\ge\bbP(J=j,\ V_{\mathrm b}\ge y)\\
&=\vartheta_j^{\mathrm b}y^{-3/4} =2^{-j-2}T_j^{-1/4}y^{-3/4}.
\end{aligned}
\]
This also holds at $y=T_j$, since the event uses $\ge$ and the truncation creates an atom at $T_j$.

\noindent Finally, $\zeta^2=Y_{\mathrm b}/4$ gives $\bbE[\zeta^2]\le1/2$, the independent symmetric sign gives $\bbE[\zeta]=0$, and $Y_{\mathrm b}\ge1$ gives $|\zeta|\ge1/2$ almost surely. For $1/2\le t\le\sqrt{T_j}/2$, we have $1\le4t^2\le T_j$, so independence also yields
\[
\bbP(\zeta\le-t) =\frac12\bbP(Y_{\mathrm b}\ge4t^2) \ge 2^{-j-3}T_j^{-1/4}(4t^2)^{-3/4}.
\]
\end{proof}

\begin{lemma}\label{lem:boundary-finite-hazard}
Fix $p>0$, an integer $M\ge1$, and a policy in $\cS_{\mathrm{fin}}$. Almost surely on \(\{Q_M<\infty\}\), one has $\Gamma_k\to0$ and $A_{T_j}=o(T_j^{1/4})=o(R_j)$.
\end{lemma}
\begin{proof}
This is a numerical argument applied to the realized predictable sequence;
we do not take expectations conditional on the future event
$\{Q_M<\infty\}$. Write $\pi_k=\bbP_k(E_k^{(M)})$ and
$b_k=\Gamma_k^2r_k^2/(16M^2)$, both $\cF_k$-measurable. Since $\zeta_k=\Xi_k\sqrt{Y_{\mathrm{b},k}}/2$, the event $E_k^{(M)}$ is $\{\Xi_k=-1,\ Y_{\mathrm{b},k}\ge1/b_k\}$, and the independence of the current sample from $\cF_k$ gives
\begin{equation}\label{eq:boundary-hazard-identity}
 \pi_k=\frac12\bbP\bigl(Y_{\mathrm{b}}\ge1/b_k\bigr),
\end{equation}
interpreted as zero when $b_k=0$. If $b_k\ge1$, then $Y_{\mathrm{b}}\ge1$ implies
$\pi_k=1/2$. Thus the set
$\mathcal J_M=\{k:b_k\ge1\}$ is finite on $\{Q_M<\infty\}$, and
\begin{equation}\label{eq:boundary-saturated-prefix}
 C_M:=\sum_{k\in\mathcal J_M}\Gamma_k<\infty.
\end{equation}
All sufficiently late indices satisfy $b_k<1$, hence
$\Gamma_k<4M/r_k\to0$. This handles arbitrarily large finite steps without imposing a cap.

\noindent Fix $j$ and write $T=T_j$. At an index $k<T$ with $k\notin\mathcal J_M$ there are two
possibilities. If $0\le b_k<T^{-1}$, then
\begin{equation}\label{eq:boundary-tiny-weight}
 \Gamma_k\le4M\,T^{-1/2}r_k^{-1}.
\end{equation}
This includes $\Gamma_k=0$. If $T^{-1}\le b_k<1$, then
$1<1/b_k\le T$, and Lemma~\ref{lem:boundary-mixture} and
\eqref{eq:boundary-hazard-identity} give
$\pi_k\ge\frac12\delta_j T^{-1/4}b_k^{3/4}$.
Solving for $\sqrt{b_k}$ and using
$\Gamma_k=4M r_k^{-1}\sqrt{b_k}$ yields
\begin{equation}\label{eq:boundary-middle-weight}
 \Gamma_k\le
 2^{8/3}M\,
 \delta_j^{-2/3}T^{1/6}r_k^{-1}\pi_k^{2/3}.
\end{equation}
The equality case $b_k=T^{-1}$ is legitimate because the tail
bound in~\eqref{eq:boundary-mixture-estimates} includes the truncation point.

\noindent Sum~\eqref{eq:boundary-tiny-weight} and
\eqref{eq:boundary-middle-weight} over their respective indices, extend
the resulting nonnegative sums to all $k<T$, and bound the indices in $\mathcal J_M$ by $C_M$. H\"older's inequality with exponents $3/2$ and $3$ gives
\[
 A_T\le C_M+4M\,T^{-1/2}\sum_{k<T}r_k^{-1}
 +2^{8/3}M\,\delta_j^{-2/3}T^{1/6}
 \Bigl(\sum_{k<T}\pi_k\Bigr)^{2/3}
 \Bigl(\sum_{k<T}r_k^{-3}\Bigr)^{1/3}.
\]
By Lemma~\ref{lem:power-sums}, $\sum_{k<T}r_k^{-1}\le C_p T^{3/4}[\log(e+T)]^{-p}$ and $\sum_{k<T}r_k^{-3}\le C_p T^{1/4}[\log(e+T)]^{-3p}$ for a deterministic constant $C_p$. Since $1/6+(1/4)/3=1/4$, enlarging $C_p$ to cover the finitely many small horizons gives
\begin{equation}\label{eq:boundary-finite-hazard-bound}
 A_{T_j}\le C_M+
 C_p M\, T_j^{1/4}[\log(e+T_j)]^{-p}
 \left(1+\delta_j^{-2/3}Q_M^{2/3}\right)
 \qquad(j\ge1).
\end{equation}
By Lemma~\ref{lem:sparse-profile}, $T_j^{1/4}/R_j\to1$ and $\log(e+T_j)\ge4\cdot2^{j^2}$, so that
$[\log(e+T_j)]^{-p}\delta_j^{-2/3}\le4^{-p}2^{-pj^2+2(j+2)/3}\to0$, as $j\to\infty$.
Dividing~\eqref{eq:boundary-finite-hazard-bound} by $T_j^{1/4}$ and using the pathwise finiteness of $Q_M$ proves $A_{T_j}=o(T_j^{1/4})=o(R_j)$.
\end{proof}


\begin{lemma}\label{lem:boundary-clock}
For every policy in $\cS_{\mathrm{fin}}$, the iterates generated by~\eqref{eq:sgd} on the instance of Section~\ref{sec:smooth-sketch} satisfy
\begin{equation}\label{eq:boundary-never-zero}
 \bbP(v_k\ne0\text{ for all }k\ge1)=1.
\end{equation}
On a path satisfying this property and $\Gamma_k\to0$, there are a
finite index $\kappa_\omega$ and a finite $C_\omega\ge0$ such that
\begin{equation}\label{eq:boundary-clock-control}
 H_2(|v_N|)\le C_\omega+3\gamma_{\mathrm{b}}A_N
 \qquad(N\ge\kappa_\omega).
\end{equation}
\end{lemma}
\begin{proof}
Projection acts separately on the coordinates, so
\[
 v_{k+1}=\proj_{[-2,2]}
 \bigl(v_k-\gamma_{\mathrm{b}}\Gamma_k\xi_k\varphi'(v_k)\bigr).
\]
If $v_k\ne0$ and $\Gamma_k=0$, a zero next coordinate is impossible.
If $v_k\ne0$ and $\Gamma_k>0$, projection can equal zero only when
its input equals zero. Since $\varphi'(v_k)\ne0$, this requires the
single value
$\xi_k=v_k/(\gamma_{\mathrm{b}}\Gamma_k\varphi'(v_k))$.
Conditionally on $\cF_k$, however, $\xi_k$ has its atomless uniform
law. Hence the conditional probability of this equality is zero.
Induction from $v_1=2$ and a countable union prove
\eqref{eq:boundary-never-zero}.

\noindent On a path with $\Gamma_k\to0$, choose a finite $\kappa_\omega$ so that
$\epsilon_k=\gamma_{\mathrm{b}}\Gamma_k\xi_k\le1/2$ for
$k\ge\kappa_\omega$. Such an index exists because $\xi_k\le3/2$.
It is only a pathwise index, and no conditional expectation is taken
after choosing it. With $w_k=|v_k|\in(0,2]$, the bound $\psi(w)\le w$ of Lemma~\ref{lem:sparse-profile}
shows that the unprojected update preserves sign and has magnitude
\[
 w_k-\epsilon_k\psi(w_k)\in[w_k/2,w_k].
\]
It stays in $[-2,2]$, so projection is inactive. Apply
Lemma~\ref{lem:clock} to $h=\psi$, $\bar b=2$, and $K_0=1$. It gives
\[
 H_2(w_{k+1})-H_2(w_k)\le2\epsilon_k
 \le3\gamma_{\mathrm{b}}\Gamma_k.
\]
Summation proves~\eqref{eq:boundary-clock-control} with
$C_\omega=H_2(w_{\kappa_\omega})<\infty$.
\end{proof}

\subsection{Proof of Theorem~\ref{thm:smooth_lb}~\ref{itm:smooth_lb_fast}}
\label{app:smooth-endpoint}

Here we provide the proof of Theorem~\ref{thm:smooth_lb}~\ref{itm:smooth_lb_fast}.
Compared with the preceding subsection, the flat coordinate is
unchanged, but the noisy coordinate now contributes a quadratic gap.
The tail exponent and the deterministic test times are adjusted
accordingly. A prescribed deterministic cap $\alpha>0$ is fixed
before the objective is selected and guarantees $\mu\Gamma_k\le1/2$, which keeps the flat coordinate positive (Lemma \ref{lem:smooth-flat-clock}).

\paragraph{Test times.} The flat coordinate is the one of Section~\ref{sec:smooth-sketch}, with the scales $\varrho_j$, $\beta_j$, $R_j$ defined there and the profile, potential, and clock of Lemma~\ref{lem:sparse-profile}. Only the test times change. Progress toward $\varrho_j$ is now tested at
\begin{equation}\label{eq:stationary-times}
 T_j^{\mathrm{s}}=\lceil R_j^3\rceil,
\end{equation}
and throughout this subsection we write $T_j=T_j^{\mathrm{s}}$. In particular, \eqref{eq:boundary-scale-growth}, which concerns $\lceil R_j^4\rceil$, is not used in this subsection; its role is played by~\eqref{eq:stationary-scale-growth}. The exponent $3$ is matched to the quadratic gap: at the test time, the critical cumulative movement $T_j^{1/3}$ is asymptotically equal to $R_j$. Precisely, $R_j^3\le T_j<R_j^3+1$ and $\log(e+T_j)\ge3\log R_j\ge3\cdot2^{j^2}$, so that, for every $p>0$,
\begin{equation}\label{eq:stationary-scale-growth}
 \begin{gathered}
 T_j^{1/3}\ge R_j,\qquad
 \frac{T_j^{1/3}}{R_j}\to1,\qquad
 \log(e+T_j)\ge3\cdot2^{j^2},\\
 \varrho_j^2\,[\log(e+T_j)]^{p}\ge3^p\,2^{pj^2-6j}\longrightarrow\infty, \ \ \text{as} \ \ j\to\infty.
 \end{gathered}
\end{equation}
As in Section~\ref{sec:smooth-sketch}, these test times do not depend on $p$.



\paragraph{Objective.}


\noindent Define
\begin{equation}\label{eq:smooth-hard-instance}
 \mu=\min\{L,(2\alpha)^{-1}\},\qquad
 F(u,v)=\mu\left(\frac{u^2}{2}+\varphi(v)\right),\qquad
 \cX=\bbR^2,\qquad x_1=(0,1/2).
\end{equation}
The curvature factor $\mu$ is fixed before the policy is run. Besides ensuring global $L$-smoothness, it guarantees $\mu\Gamma_k\le1/2$ for every policy in $\cS_\alpha$, which is used in Lemma~\ref{lem:smooth-flat-clock}.

\begin{lemma}\label{lem:smooth-hard-regularity}
The objective \eqref{eq:smooth-hard-instance} is convex and $C^2$, has
globally $L$-Lipschitz gradient, has unique minimizer $(0,0)$ with
$F_\star=0$, and is not strongly convex.
\end{lemma}
\begin{proof}
Lemma~\ref{lem:sparse-profile} gives $0\le\varphi''\le1$ and
$\varphi''(0)=0$. Hence
\[
 \nabla^2F(u,v)=\mu\begin{pmatrix}1&0\\0&\varphi''(v)\end{pmatrix},
 \qquad 0\preceq\nabla^2F(u,v)\preceq\mu I\preceq LI.
\]
Both summands are nonnegative and vanish only at zero, proving uniqueness.
The Hessian at zero has a zero eigenvalue, precluding a positive global
strong-convexity modulus. Its gradient at zero vanishes.
\end{proof}

\paragraph{First-Order Oracle.} Here we define a fixed finite-variance noise law; Lemma~\ref{lem:smooth-flat-clock} below verifies Assumption~\ref{ass:abc_condition}. Let
\begin{equation}\label{eq:smooth-mixture-weights}
 \delta_j=2^{-j-2},\qquad \vartheta_j=\delta_jT_j^{-1/3}.
\end{equation}
Their sum is at most $\sum_j\delta_j=1/4$.
Let $V\ge1$ have survival function $\bbP(V\ge y)=y^{-2/3}$ for $y\ge1$.
Define the distribution of $Y$ as follows: with mass $\vartheta_j$, use
$\min\{V,T_j\}$, and put the remaining mass $1-\sum_j\vartheta_j$ at
$1$. The law of $Y$ depends neither on $k$ nor on the policy.

\noindent The mixing weights have been chosen so that each truncation contributes
only a summable amount to the variance, while still supplying a tail lower
bound at its own prescribed scale. 

\begin{lemma}\label{lem:smooth-mixture}
For an independent Rademacher sign $\Xi$, let $\zeta=\Xi\sqrt Y$.
Then $\bbE[\zeta]=0$ and
\begin{equation}\label{eq:smooth-mixture-moment}
 1\le\sigma_\zeta^2:=\bbE[\zeta^2]=\bbE[Y]
 =1+3\sum_j\delta_j-3\sum_j\vartheta_j\le\frac74.
\end{equation}
Moreover, for every $j$,
\begin{equation}\label{eq:smooth-mixture-tail}
 \bbP(Y\ge y)\ge\delta_jT_j^{-1/3}y^{-2/3}
 \qquad(1\le y\le T_j),
\end{equation}
and $Y\ge1$ almost surely.
\end{lemma}
\begin{proof}
The layer-cake formula gives
$\bbE[\min\{V,T_j\}]=1+\int_1^{T_j}y^{-2/3}\,dy=3T_j^{1/3}-2$.
Multiplying by the component masses and adding the atom at one proves
\eqref{eq:smooth-mixture-moment}; interchange of the nonnegative series
and integral is valid by monotone convergence. Finite second moment
implies finite first absolute moment of $\zeta$, so symmetry gives its
zero mean. Keeping only the $j$th mixture component proves
\eqref{eq:smooth-mixture-tail}. At $y=T_j$ that component includes the
atom created by truncation, so the endpoint inequality is also valid.
\end{proof}

\noindent Take i.i.d. copies $\zeta_k$, independent of the algorithmic seed $U$, which is $\cF_1$-measurable (Appendix~\ref{app:details_predictable_ss}),
and set $\cF_k=\sigma(U,\zeta_1,\ldots,\zeta_{k-1})$, completed as
usual; any policy in $\cS_\alpha$ is predictable with respect to this filtration, since all information available before the $k$th oracle call is generated by $U$ and $\zeta_1,\dots,\zeta_{k-1}$. The fixed additive oracle is
\begin{equation}\label{eq:smooth-hard-oracle}
 \widehat{\nabla}_k F(X_k)=(\mu u_k+\zeta_k,\,\mu\varphi'(v_k)).
\end{equation}
It is the gradient of the random convex sample loss $F(u,v)+\zeta_k u$.
The recursions are
\begin{equation}\label{eq:smooth-hard-dynamics}
 u_{k+1}=(1-\mu\Gamma_k)u_k-\Gamma_k\zeta_k,
 \qquad v_{k+1}=v_k-\mu\Gamma_k\psi(v_k).
\end{equation}
The second identity is valid because the next lemma keeps $v_k>0$.

\noindent Next, we check that the cap keeps the flat coordinate positive and that the oracle satisfies Assumption~\ref{ass:abc_condition}.

\begin{lemma}\label{lem:smooth-flat-clock}
For every policy in $\cS_\alpha$, it holds that
\begin{equation}\label{eq:smooth-flat-clock}
 0<v_k\le1/2,\qquad v_k/2\le v_{k+1}\le v_k,\qquad
 H_{1/2}(v_N)\le2\mu A_N.
\end{equation}
The iterates have finite second moments at each finite time, and
\begin{equation}\label{eq:smooth-hard-moments}
 \bbE_k[\widehat{\nabla}_k F(X_k)]=\nabla F(X_k),\qquad
 \bbE_k[\|{\widehat{\nabla}_k F(X_k)}\|^2]=\norm{\nabla F(X_k)}^2+\sigma_\zeta^2.
\end{equation}
In particular, Assumption~\ref{ass:abc_condition} holds with
$a=0,b=1,c=\sigma_\zeta^2$.
\end{lemma}
\begin{proof}
Since $0\le\mu\Gamma_k\le1/2$, Lemma~\ref{lem:clock} with $h=\psi$, $\bar b=1/2$, and $K_0=1$ gives positivity and monotonicity, starting
from $v_1=1/2$. Its clock inequality telescopes from $H_{1/2}(v_1)=0$
to give the final bound in \eqref{eq:smooth-flat-clock}.
Also $|1-\mu\Gamma_k|\le1$, so
$|u_{k+1}|\le|u_k|+\alpha|\zeta_k|$. The inequality
$u_{k+1}^2\le2u_k^2+2\alpha^2\zeta_k^2$ proves finite second moments
by induction. Conditional independence and centering of $\zeta_k$ now
justify both identities in \eqref{eq:smooth-hard-moments}.
\end{proof}


\noindent By Lemma~\ref{lem:smooth-flat-clock}, slow progress of the flat coordinate reduces to a bound on the cumulative stepsize, which Lemma~\ref{lem:finite-hazard} derives from a finite shock hazard; Lemma~\ref{lem:smooth-spike} handles the complementary case.

\noindent Fix $p>0$ and put $r_k=k^{1/3}[\log(e+k)]^{p}$, the reciprocal of the rate in Theorem~\ref{thm:smooth_lb}~\ref{itm:smooth_lb_fast}. For every integer $M\ge1$,
define
\begin{equation}\label{eq:smooth-hazard-def}
 E_k^{(M)}=\{\Gamma_k|\zeta_k|\ge M/\sqrt{r_k}\},\quad
 \pi_k^{(M)}=\bbP_k(E_k^{(M)}),\quad
 Q_M=\sum_{k\ge1}\pi_k^{(M)}.
\end{equation}
The event is in $\cF_{k+1}$. Letting $b_k=\Gamma_k^2r_k/M^2$, the
independence of the new noise gives
\begin{equation}\label{eq:smooth-hazard-tail}
 \pi_k^{(M)}=\bbP(Y\ge1/b_k),
\end{equation}
where this denotes evaluation of the deterministic tail at a predictable
random threshold and equals zero if $b_k=0$.

\noindent The next lemma converts a finite shock hazard into a bound on the cumulative stepsize at the scale $T_j^{1/3}$; the cap $\alpha$ controls thresholds beyond the truncation point.

\begin{lemma}\label{lem:finite-hazard}
For every fixed $p>0$ and $M\ge1$, there is a deterministic $j_0$ such
that, for every policy in $\cS_\alpha$ and every $j\ge j_0$,
\begin{equation}\label{eq:finite-hazard-quantitative}
 A_{T_j}\le C_p M T_j^{1/3}[\log(e+T_j)]^{-p/2}
 \left[1+\delta_j^{-3/4}
               \left(\sum_{k<T_j}\pi_k^{(M)}\right)^{3/4}\right].
\end{equation}
Here $C_p$ is deterministic and finite; $j_0$ may also depend on
$\alpha$ and $M$. Consequently,
\begin{equation}\label{eq:finite-hazard-implication}
 Q_M<\infty\quad\Longrightarrow\quad
 A_{T_j}=o(T_j^{1/3})=o(R_j).
\end{equation}
\end{lemma}
\begin{proof}
For this proof let $N=T_j$. For every $k<N$, the deterministic cap and
monotonicity of $r_k$ give
\[
 b_k\le(\alpha/M)^2N^{1/3}[\log(e+N)]^{p}.
\]
Thus
\begin{equation}\label{eq:saturated-tail-check}
 \delta_jN^{-1/3}b_k^{2/3}
 \le\delta_j(\alpha/M)^{4/3}N^{-1/9}[\log(e+N)]^{2p/3}
 \longrightarrow0, \ \ \text{as} \ \ N\to\infty
\end{equation}
uniformly over all $k<N$ and all policies. Choose $j_0$ so that this
upper bound is at most one for $j\ge j_0$.
If $N^{-1}\le b_k\le1$, the tail bound
\eqref{eq:smooth-mixture-tail} applies at $1/b_k\in[1,N]$. If $b_k>1$,
$Y\ge1$ makes $\pi_k^{(M)}=1$, and \eqref{eq:saturated-tail-check}
gives the same conclusion. Hence in both regimes
\begin{equation}\label{eq:hazard-lower-at-scale}
 b_k\ge N^{-1}\quad\Longrightarrow\quad
 \pi_k^{(M)}\ge\delta_jN^{-1/3}b_k^{2/3}.
\end{equation}
For $b_k<N^{-1}$, one has $\Gamma_k\le M N^{-1/2}r_k^{-1/2}$.
For $b_k\ge N^{-1}$, solving \eqref{eq:hazard-lower-at-scale} for
$\sqrt{b_k}$ gives
$\Gamma_k\le M\delta_j^{-3/4}N^{1/4}r_k^{-1/2}(\pi_k^{(M)})^{3/4}$.
Combining these alternatives, including $\Gamma_k=0$, yields
\[
 \Gamma_k\le M N^{-1/2}r_k^{-1/2}
 +M\delta_j^{-3/4}N^{1/4}r_k^{-1/2}(\pi_k^{(M)})^{3/4}.
\]
Sum over $k<N$ and apply H\"older's inequality with exponents $4/3$ and
$4$ to the second sum:
\begin{align*}
 A_N&\le M N^{-1/2}\sum_{k<N}r_k^{-1/2}+M\delta_j^{-3/4}N^{1/4}\left(\sum_{k<N}\pi_k^{(M)}\right)^{3/4}\left(\sum_{k<N}r_k^{-2}\right)^{1/4}.
\end{align*}
Lemma~\ref{lem:power-sums} gives
\[
 \sum_{k<N}r_k^{-1/2}\le C_pN^{5/6}[\log(e+N)]^{-p/2},\qquad
 \sum_{k<N}r_k^{-2}\le C_pN^{1/3}[\log(e+N)]^{-2p}.
\]
Substitution proves \eqref{eq:finite-hazard-quantitative}.
If $Q_M<\infty$, replace its finite partial sum by the fixed
finite pathwise number $Q_M$. Since $\log(e+T_j)\ge3\cdot2^{j^2}$ by~\eqref{eq:stationary-scale-growth},
\[
 [\log(e+T_j)]^{-p/2}\delta_j^{-3/4}
 \le C_p2^{-pj^2/2+3j/4}\longrightarrow0, \ \ \text{as} \ \ j\to\infty.
\]
Both terms in the bracketed estimate after division by $T_j^{1/3}$ tend
to zero. Finally, $T_j^{1/3}/R_j\to1$ by~\eqref{eq:stationary-scale-growth}. This proves
\eqref{eq:finite-hazard-implication} pathwise; it does not assume a
deterministic bound on $Q_M$.
\end{proof}

\noindent The preceding estimate concerns probabilities and movement, rather than
the objective at an actual spike time. A cancellation in the quadratic
coordinate might hide a large innovation at the next iterate; the following
two-iterate argument shows that the current iterate must then already have
a correspondingly large gap.

\begin{lemma}\label{lem:smooth-spike}
Every occurrence of $E_k^{(M)}$ implies
\begin{equation}\label{eq:smooth-spike-transfer}
 \max\{r_k\Delta(X_k),r_{k+1}\Delta(X_{k+1})\}\ge\frac{\mu M^2}{8}.
\end{equation}
\end{lemma}
\begin{proof}
Let $d=\Gamma_k|\zeta_k|\ge M/\sqrt{r_k}$. If $|u_k|\ge d/2$, then
$r_k\Delta(X_k)\ge\mu r_ku_k^2/2\ge\mu M^2/8$.
Otherwise, \eqref{eq:smooth-hard-dynamics} and $|1-\mu\Gamma_k|\le1$
give $|u_{k+1}|\ge d-|u_k|>d/2$.
As $r_{k+1}\ge r_k$, the same bound holds at $k+1$. No favorable sign
of the innovation is required.
\end{proof}


\noindent The preceding estimates complete the two sides of the argument.
We partition the analysis according to cumulative movement at
the prescribed scales: small movement leaves flat-coordinate error,
whereas the complementary regime forces divergent spike hazard
at every amplitude. 

\begin{proof}[Theorem~\ref{thm:smooth_lb}~\ref{itm:smooth_lb_fast}]
Lemmas~\ref{lem:smooth-hard-regularity}, \ref{lem:smooth-mixture}, and
\ref{lem:smooth-flat-clock} verify all instance and oracle assertions.
Fix a policy in $\cS_\alpha$ and first fix $p>0$. Define the measurable
partition
\[
 \Omega_{\mathrm{small}}=\{\liminf_j A_{T_j}/R_j=0\},\qquad
 \Omega_{\mathrm{large}}=\{\liminf_j A_{T_j}/R_j>0\}.
\]

\noindent On $\Omega_{\mathrm{small}}$ choose a pathwise subsequence $j_m\to\infty$
with $A_{T_{j_m}}/R_{j_m}\to0$. The clock bound gives
$H_{1/2}(v_{T_{j_m}})=o(R_{j_m})$. Since
$H_{1/2}(\varrho_{j_m})\ge R_{j_m}$ and $H_{1/2}$ is strictly
decreasing, eventually $v_{T_{j_m}}>\varrho_{j_m}$.
Therefore
\begin{align*}
 r_{T_{j_m}}\Delta(X_{T_{j_m}})
 &\ge\frac{\mu}{2}T_{j_m}^{1/3}[\log(e+T_{j_m})]^{p}
                     \varrho_{j_m}\beta_{j_m}\\
 &\ge\frac{\mu}{2}\varrho_{j_m}^2[\log(e+T_{j_m})]^{p}
 \longrightarrow\infty, \ \ \text{as} \ \ j_m\to\infty
\end{align*}
by $T_{j_m}\ge R_{j_m}^3$, the plateau bound of Lemma~\ref{lem:sparse-profile}, and~\eqref{eq:stationary-scale-growth}. This part is purely pathwise.

\noindent On $\Omega_{\mathrm{large}}$, Lemma~\ref{lem:finite-hazard} implies
$Q_M=\infty$ for every integer $M\ge1$, since finite hazard
would make the liminf zero. Lemma~\ref{lem:cbc} shows, for each such $M$,
that $E_k^{(M)}$ occurs infinitely often almost surely on this event.
Intersecting these countably many probability-one statements, Lemma~\ref{lem:smooth-spike} yields, for every $M$, infinitely many indices $k$ with $\max\{r_k\Delta(X_k),r_{k+1}\Delta(X_{k+1})\}\ge\mu M^2/8$; hence $\limsup_k r_k\Delta(X_k)=\infty$. The two events cover every path, proving the result for fixed
$p>0$.

\noindent Apply this proof to $p=1/m$, $m\in\bbN$, and intersect the resulting
probability-one events. For every real $p>0$, choose $m$ with $1/m<p$.
Because $\Delta(X_k)\ge0$ and $[\log(e+k)]^{p-1/m}\ge1$, the conclusion at
$1/m$ implies the conclusion at $p$. Thus all positive logarithmic
exponents hold on one event for the fixed policy. The instance depends neither on $p$ nor on the policy. This completes the proof.
\end{proof}



\newpage
\section{Proof of Theorem~\ref{thm:smooth_ub}}\label{app:upper-smooth}

We turn from lower-bound construction to the {upper-bound proof}.
We first establish moment and weighted-descent estimates for the
projected iterates. The argument follows the standard almost-supermartingale framework
\citep{Robbins1971,CombettesPesquet2015}, with explicit growth bounds
under the abc condition. These bounds allow a nonzero gradient
at the constrained minimizer.

\begin{lemma}\label{lem:smooth-energy}
Suppose Assumptions~\ref{ass:convexity}, \ref{ass:domain},
\ref{ass:smoothness}, and~\ref{ass:abc_condition} hold. Let $\Gamma_k$
be deterministic and positive with $S_2:=\sum_k\Gamma_k^2<\infty$.
Fix $x_\star\in\cX_\star$, and put $g_\star=\norm{\nabla F(x_\star)}$.
Define the deterministic constants
\[
 d_0=(a/2+2b)g_\star^2+c,\qquad
 d_1=a(1+L)/2+2bL^2,\qquad
 M_V=(\|x_1-x_\star\|^2+d_0S_2)e^{d_1S_2}.
\]
Then
\begin{align}
 &\bbE_k[\|{\widehat{\nabla}_k F(X_k)}\|^2]\le d_0+d_1\norm{X_k-x_\star}^2,\label{eq:smooth-Q-growth}\\
 &\bbE_k\left[\norm{X_{k+1}-x_\star}^2\right]
 \le(1+d_1\Gamma_k^2)\norm{X_k-x_\star}^2-2\Gamma_k\Delta(X_k)+d_0\Gamma_k^2,
 \label{eq:smooth-distance-general}\\
 &\sup_k\bbE[\norm{X_k-x_\star}^2]\le M_V,\quad
 2\sum_k\Gamma_k\bbE[\Delta(X_k)]
 \le \norm{X_1-x_\star}^2+(d_1M_V+d_0)S_2<\infty.\label{eq:smooth-energy-expected}
\end{align}
Moreover, $\norm{X_k-x_\star}^2$ converges almost surely to a finite limit, and
\begin{equation}\label{eq:smooth-uniform-moments}
 \begin{gathered}
 \sup_k\bbE\!\left[\bbE_k[\|{\widehat{\nabla}_k F(X_k)}\|^2]\right]<\infty,\qquad
 \sup_k\bbE[\norm{\nabla F(X_k)}^2]<\infty,\\
 \sup_k\norm{\nabla F(X_k)}^2<\infty\quad\as
 \end{gathered}
\end{equation}
\end{lemma}
\begin{proof}
The smoothness inequality at $x_\star$, followed by Young's inequality,
gives, for $x\in\cX$,
\[
 0\le\Delta(x)\le g_\star\norm{x-x_\star}
                      +\frac L2\norm{x-x_\star}^2
 \le\frac{g_\star^2}{2}+\frac{1+L}{2}\norm{x-x_\star}^2.
\]
The gradient Lipschitz property also gives
$\norm{\nabla F(x)}^2\le2g_\star^2+2L^2\norm{x-x_\star}^2$.
Substituting these two bounds into the abc inequality proves
\eqref{eq:smooth-Q-growth}. Projection nonexpansiveness and convexity
imply
\begin{align*}
 \bbE_k\left[\norm{X_{k+1}-x_\star}^2\right]
 &\le \norm{X_k-x_\star}^2-2\Gamma_k\ip{\nabla F(X_k)}{X_k-x_\star}+\Gamma_k^2\bbE_k[\|{\widehat{\nabla}_k F(X_k)}\|^2]\\
 &\le \norm{X_k-x_\star}^2-2\Gamma_k\Delta(X_k)+\Gamma_k^2(d_0+d_1\norm{X_k-x_\star}^2),
\end{align*}
which is \eqref{eq:smooth-distance-general}.
All required moments are finite by induction. Indeed, $\bbE\bigl[\norm{X_k-x_\star}^2\bigr]<\infty$
first bounds the expected gap and the expected conditional second moment
through the two growth inequalities above. The distance recursion then
gives $\bbE\bigl[\norm{X_{k+1}-x_\star}^2\bigr]<\infty$.
Dropping descent and iterating expectations yields
\[
 \bbE\bigl[\norm{X_k-x_\star}^2\bigr]\le\left(\norm{X_1-x_\star}^2+d_0\sum_{i<k}\Gamma_i^2\right)
                 \prod_{i<k}(1+d_1\Gamma_i^2)\le M_V.
\]
Summing the expected distance inequality and dropping the final
nonnegative distance gives the weighted energy bound in
\eqref{eq:smooth-energy-expected}. Theorem~\ref{lem:rs}, with
$G_k=\norm{X_k-x_\star}^2$, $a_k=d_1\Gamma_k^2$, $D_k=2\Gamma_k\Delta(X_k)$, and
$B_k=d_0\Gamma_k^2$, proves finite almost-sure convergence of $\norm{X_k-x_\star}^2$.
Its entire sequence is consequently bounded on that probability-one
event. The gradient growth bound and \eqref{eq:smooth-Q-growth} now
prove \eqref{eq:smooth-uniform-moments}. The pathwise bound is used only
after the probabilistic estimates have been completed.
\end{proof}


\noindent Both proofs below start from this energy estimate. They differ in the stepsize schedule and in how the smoothing bias of the Moreau envelope is controlled: without a stationary minimizer, the bias is additive, so the smoothing parameter must decrease; with a stationary minimizer, the bias is proportional to the gap, so a fixed smoothing parameter suffices.

\noindent Before giving the formal proofs, recall that $\Phi=\operatorname{cl}(F+\iota_\cX)$, the function whose epigraph is the closure of that of $F+\iota_\cX$. Given a \emph{smoothing parameter} $\lambda>0$, the Moreau envelope of $\Phi$ \citep{BauschkeCombettes2017} is defined as

\[
\Phi_\lambda(x) \coloneqq \min_{y\in\cH}\Bigl\{\Phi(y)+\frac{\norm{y-x}^2}{2\lambda}\Bigr\},\qquad x\in\cH,
\]
with the properties of Lemmas~\ref{lem:closed-envelope}-\ref{lem:envelope-step}.

\subsection{Proof of Theorem~\ref{thm:smooth_ub}~\ref{itm:smooth_ub_slow}}\label{app:upper-smooth-slow}
Fix $\eta>0$ and set $p=1/2+\eta$. Choose a deterministic amplitude $\gamma_0>0$ and a deterministic $t_0\ge2$ with $\log(1+t_0)>4p$, and run~\eqref{eq:sgd} with the stepsizes
\begin{equation}\label{eq:stepsize-slow}
 \Gamma_k=\gamma_0(k+t_0)^{-3/4}[\log(k+t_0)]^{-p},
\end{equation}
which are those of Lemma~\ref{lem:quarter-schedule}. They are deterministic, positive, and decreasing in $k$, hence belong to $\cS_\alpha$ with $\alpha=\Gamma_1$; any smaller prescribed positive cap can be met by increasing $t_0$. The analysis also uses the smoothing parameters $\lambda_k$ and the normalization $r_k$ of Lemma~\ref{lem:quarter-schedule}, with an arbitrary $\lambda_0>0$; here $r_k$ is the reciprocal of the rate in Theorem~\ref{thm:smooth_ub}~\ref{itm:smooth_ub_slow}. The moment estimates of Lemma~\ref{lem:smooth-energy} justify the conditional envelope recursion, the schedule makes its normalized perturbations summable in expectation, and the pathwise gradient bound removes the smoothing bias.

\begin{proof}[Theorem~\ref{thm:smooth_ub}~\ref{itm:smooth_ub_slow}]
Since $\sum_k\Gamma_k^2<\infty$ by Lemma~\ref{lem:quarter-schedule}, Lemma~\ref{lem:smooth-energy} applies. All envelope variables are
measurable by Lemma~\ref{lem:envelope} and satisfy
$0\le \Phi_{\lambda_k}(X_k)-F_\star\le\Delta(X_k)$ by Lemmas~\ref{lem:closed-envelope} and~\ref{lem:envelope-estimates}, which ensures their integrability by
Lemma~\ref{lem:smooth-energy}.

\noindent Lemma~\ref{lem:envelope-step} gives
\[
 \bbE_k\bigl[\Phi_{\lambda_{k+1}}(X_{k+1})-F_\star\bigr]\le \Phi_{\lambda_k}(X_k)-F_\star+\frac{\Gamma_k^2}{2\lambda_{k+1}}\bbE_k[\|{\widehat{\nabla}_k F(X_k)}\|^2]
       +\frac{\lambda_k-\lambda_{k+1}}2\norm{\nabla F(X_k)}^2.
\]
Let $M_Q=\sup_k\bbE\!\left[\bbE_k[\|{\widehat{\nabla}_k F(X_k)}\|^2]\right]$ and
$M_g=\sup_k\bbE[\norm{\nabla F(X_k)}^2]$, both finite by
\eqref{eq:smooth-uniform-moments}. Since all scalar coefficients are
deterministic,
\begin{align*}
 \sum_k r_{k+1}\bbE\left[\frac{\Gamma_k^2}{2\lambda_{k+1}}\bbE_k[\|{\widehat{\nabla}_k F(X_k)}\|^2]
       +\frac{\lambda_k-\lambda_{k+1}}2\norm{\nabla F(X_k)}^2\right]\\
 \le\frac12\sum_k r_{k+1}
 \left[M_Q\frac{\Gamma_k^2}{\lambda_{k+1}}
          +M_g(\lambda_k-\lambda_{k+1})\right]<\infty.
\end{align*}
by~\eqref{eq:quarter-sums}. Also $\sum_k\Gamma_k\bbE\bigl[\Phi_{\lambda_k}(X_k)-F_\star\bigr]\le\sum_k\Gamma_k\bbE[\Delta(X_k)]<\infty$ by~\eqref{eq:smooth-energy-expected}, and $r_k$ is nondecreasing with $r_{k+1}-r_k\le\Gamma_k/(4\gamma_0)$ and $\sum_k\Gamma_k/r_k=\infty$ by Lemma~\ref{lem:quarter-schedule}.
All hypotheses of Lemma~\ref{lem:scaled} therefore hold with $a_k=0$ and
$\theta_k=\Gamma_k$, so $r_k\bigl(\Phi_{\lambda_k}(X_k)-F_\star\bigr)\to0$ almost surely.

\noindent The envelope gap bound~\eqref{eq:envelope-gap} in Lemma~\ref{lem:envelope-estimates} yields
\[
 0\le r_k\Delta(X_k)-r_k\bigl(\Phi_{\lambda_k}(X_k)-F_\star\bigr)
 \le\frac{r_k\lambda_k}{2}\norm{\nabla F(X_k)}^2.
\]
The first factor tends to zero by \eqref{eq:quarter-errors}, and the
second is pathwise bounded by \eqref{eq:smooth-uniform-moments}.
Consequently $r_k\Delta(X_k)\to0$ almost surely, proving the stated result.
\end{proof}

\subsection{Proof of Theorem~\ref{thm:smooth_ub}~\ref{itm:smooth_ub_fast}}\label{app:upper-smooth-fast}
Suppose, in addition, that Assumption~\ref{ass:stationary_minimizer} holds. Fix $\eta>0$ and set $p=1/3+\eta$. Choose a deterministic amplitude $\gamma_0>0$, put $K=a+2bL$, and select a deterministic $t_0\ge2$ large enough that
\[
 \log(1+t_0)>3p,\qquad
 \frac{K\gamma_0}{(1+t_0)^{2/3}[\log(1+t_0)]^p}\le1.
\]
The second inequality is automatic when $K=0$. Both requirements hold for all sufficiently large $t_0$,
because both $\log(1+t_0)$ and $(1+t_0)^{2/3}[\log(1+t_0)]^p$ tend to infinity as $t_0\to\infty$.
Run~\eqref{eq:sgd} with the stepsizes
\begin{equation}\label{eq:stepsize-fast}
 \Gamma_k=\frac{\gamma_0}{(k+t_0)^{2/3}[\log(k+t_0)]^p}.
\end{equation}
They are deterministic, positive, and decreasing in $k$, hence belong to $\cS_\alpha$ with $\alpha=\Gamma_1$; increasing $t_0$ meets any prescribed positive cap. With a stationary minimizer, the standard gradient self-bound from smooth convex analysis \citep{BauschkeCombettes2017} controls the envelope approximation error by a fixed fraction of the objective gap. This permits a fixed smoothing parameter and the faster one-third schedule. We first establish the comparison and its conditional recursion.

\begin{lemma}\label{lem:stationary-envelope}
Under Assumptions~\ref{ass:convexity}, \ref{ass:domain},
\ref{ass:smoothness}, \ref{ass:abc_condition}, and~\ref{ass:stationary_minimizer},
\begin{equation}\label{eq:gradient-self-bound}
 \norm{\nabla F(x)}^2\le2L(F(x)-F_\star)\qquad(x\in\cH).
\end{equation}
For $0<\lambda<1/L$, $\Phi=F+\iota_\cX$, we have
\begin{equation}\label{eq:stationary-envelope-equivalence}
 (1-L\lambda)\Delta(X_k)\le \Phi_\lambda(X_k)-F_\star\le\Delta(X_k).
\end{equation}
Set $K=a+2bL$. Then
\begin{equation}\label{eq:stationary-envelope-recursion}
 \bbE_k[\Phi_\lambda(X_{k+1})-F_\star]\le(1+c_1\Gamma_k^2)\bigl(\Phi_\lambda(X_k)-F_\star\bigr)+c_2\Gamma_k^2,
 \quad c_1=\frac{K}{2\lambda(1-L\lambda)},\quad c_2=\frac c{2\lambda},
\end{equation}
whenever the terms are integrable.
\end{lemma}
\begin{proof}
The supporting inequality at a stationary $x_\star$ implies
$F(x)\ge F_\star$ on the whole Hilbert space. The descent inequality at
$x-L^{-1}\nabla F(x)$ gives
\[
 F_\star\le F(x-L^{-1}\nabla F(x))
 \le F(x)-\frac1{2L}\norm{\nabla F(x)}^2,
\]
proving \eqref{eq:gradient-self-bound}. Lemma~\ref{lem:envelope-estimates} then
gives $0\le\Delta(X_k)-\bigl(\Phi_\lambda(X_k)-F_\star\bigr)\le\lambda\norm{\nabla F(X_k)}^2/2
\le L\lambda\Delta(X_k)$, which is
\eqref{eq:stationary-envelope-equivalence}.
Furthermore $\bbE_k[\|{\widehat{\nabla}_k F(X_k)}\|^2]\le K\Delta(X_k)+c$. Apply
\eqref{eq:fixed-envelope-step}, discard its nonpositive term, and use
$\Delta(X_k)\le \bigl(\Phi_\lambda(X_k)-F_\star\bigr)/(1-L\lambda)$ to obtain
\eqref{eq:stationary-envelope-recursion}.
\end{proof}

\begin{proof}[Theorem~\ref{thm:smooth_ub}~\ref{itm:smooth_ub_fast}]
Set $r_k=(k+t_0)^{1/3}[\log(k+t_0)]^{-p}$, the reciprocal of the rate in Theorem~\ref{thm:smooth_ub}~\ref{itm:smooth_ub_fast}, and fix
any $\lambda\in(0,1/L)$. The series $\sum_k\Gamma_k^2$ converges, so
Lemma~\ref{lem:smooth-energy} supplies finite moments at each
step. Select a stationary $x_\star$,
projection and \eqref{eq:gradient-self-bound} then give
\begin{align*}
 \bbE_k\bigl[\norm{X_{k+1}-x_\star}^2\bigr]
 &\le \norm{X_k-x_\star}^2-2\Gamma_k\Delta(X_k)+\Gamma_k^2(K\Delta(X_k)+c)\\
 &\le \norm{X_k-x_\star}^2-\Gamma_k\Delta(X_k)+c\Gamma_k^2,
\end{align*}
where $K\Gamma_k\le K\Gamma_1\le1$ follows from the schedule.
Taking expectations and summing gives
\begin{equation}\label{eq:stationary-energy}
 \sum_k\Gamma_k\bbE[\Delta(X_k)]\le \norm{X_1-x_\star}^2+c\sum_k\Gamma_k^2<\infty.
\end{equation}
By Lemma~\ref{lem:envelope}, $\Phi_\lambda(X_k)-F_\star$ is $\cF_k$-measurable, and $0\le \Phi_\lambda(X_k)-F_\star\le\Delta(X_k)$ by~\eqref{eq:stationary-envelope-equivalence}, so $\Phi_\lambda(X_k)-F_\star$ is integrable by Lemma~\ref{lem:smooth-energy}. Equation~\eqref{eq:stationary-envelope-recursion} then gives the recursion
required by Lemma~\ref{lem:scaled}, with
$a_k=c_1\Gamma_k^2$, $B_k=c_2\Gamma_k^2$, and $\theta_k=\Gamma_k$.
For $R(x)=x^{1/3}(\log x)^{-p}$,
\[
 R'(x)=x^{-2/3}(\log x)^{-p}
          \left(\frac13-\frac p{\log x}\right).
\]
Thus $r_k$ increases because $\log(1+t_0)>3p$, and integration over $[k+t_0,k+1+t_0]$ gives
$0\le r_{k+1}-r_k\le\Gamma_k/(3\gamma_0)$.
The bounded ratios at successive indices also give
\[
 \sum_k r_{k+1}\Gamma_k^2
 \le C\gamma_0^2\sum_k\frac1{(k+t_0)[\log(k+t_0)]^{3p}}<\infty,
 \qquad \frac{\Gamma_k}{r_k}=\frac{\gamma_0}{k+t_0}.
\]
The series is finite because $3p>1$; the last ratio has divergent sum.
Since $0\le \Phi_\lambda(X_k)-F_\star\le\Delta(X_k)$, \eqref{eq:stationary-energy} supplies the
remaining hypothesis. Lemma~\ref{lem:scaled} proves $r_k\bigl(\Phi_\lambda(X_k)-F_\star\bigr)\to0$.
The positive constant $1-L\lambda$ in
\eqref{eq:stationary-envelope-equivalence} then gives
$r_k\Delta(X_k)\to0$ almost surely, which is the stated rate in
Theorem~\ref{thm:smooth_ub}~\ref{itm:smooth_ub_fast}, because
$p=1/3+\eta$.
\end{proof}


\newpage
\section{Proof of Theorem~\ref{thm:nonsmooth_lb}}\label{app:nonsmooth-endpoint}

This appendix proves Theorem~\ref{thm:nonsmooth_lb}. The construction uses the flat coordinate of Section~\ref{sec:smooth-sketch} and replaces the smooth first coordinate by an absolute value. Unlike the smooth constructions, the oracle is calibrated after the predictable stepsize is chosen. We first verify that this protocol is admissible, and then treat separately the two possible behaviors of the weighted square sum $\sum_kk^{1/2}\Gamma_k^2$; this separation is used only in the analysis.

\paragraph{Objective.} The flat coordinate is the one of Section~\ref{sec:smooth-sketch}: the scales $\varrho_j$, $\beta_j$, $R_j$ defined there, together with the profile $\psi$, the potential $\varphi$, and the clock $H_{\bar b}$ of Lemma~\ref{lem:sparse-profile}. Since the absolute-value coordinate also produces a linear gap, the test times are those of Theorem~\ref{thm:smooth_lb}~\ref{itm:smooth_lb_slow},
\begin{equation}\label{eq:nonsmooth-times}
 T_j^{\mathrm{n}}=\lceil R_j^4\rceil,
\end{equation}
so that the critical cumulative movement $(T_j^{\mathrm{n}})^{1/4}$ is again asymptotically equal to $R_j$. In particular, \eqref{eq:boundary-scale-growth} in Lemma~\ref{lem:sparse-profile} gives, for every $p>0$,
\begin{equation}\label{eq:nonsmooth-scale-growth}
 (T_j^{\mathrm{n}})^{1/4}\ge R_j,\qquad
 \lim_{j\to\infty}\frac{(T_j^{\mathrm{n}})^{1/4}}{R_j}=1,\qquad
 \varrho_j^2\,[\log(e+T_j^{\mathrm{n}})]^{p}\longrightarrow\infty, \ \ \text{as} \ \ j\to\infty.
\end{equation}



\noindent Let
\begin{equation}\label{eq:nonsmooth-hard-instance}
 \gamma_{\mathrm{n}}=\frac12,\qquad
 F(u,v)=\gamma_{\mathrm{n}}\bigl(|u|+\varphi(v)\bigr),\qquad
 \cX=[-1,1]\times[-2,2],\qquad x_1=(0,2).
\end{equation}
This function is convex and nonsmooth, with unique minimizer $(0,0)$
and $F_\star=0$. Its subgradients have norm at most
$\tfrac12\sqrt{1+\beta_1^2}<1$, by Lemma~\ref{lem:sparse-profile}.
Thus it is globally $1$-Lipschitz in the Euclidean norm. Fix the
subgradient selection $\sgn(0)=0$.

\paragraph{First-Order Oracle.} After $\Gamma_k$ has been selected, form the predictable quantities
\begin{equation}\label{eq:nonsmooth-protocol-scales}
 \begin{gathered}
 a_k^{\mathrm{n}}=k^{1/2}\Gamma_k^2,\quad
 b_k=\min\{1/2,a_k^{\mathrm{n}}\},\quad
 d_k=\sqrt{1+\sum_{i=1}^k b_i},\quad c_k=\sqrt{d_k},\\
 h_k=\min\{1/2,c_kk^{-1/4}\},\qquad
 \pi_k=\min\{1/4,\Gamma_k^2/(2h_k^2)\}.
 \end{gathered}
\end{equation}
   Conditionally on $\cF_k$, draw independently $J_k\in\{1,-1,0\}$ with probabilities $\pi_k/2,\pi_k/2,1-\pi_k$ and $\xi_k\sim\mathrm{Unif}[1/2,3/2]$; both are $\cF_{k+1}$-measurable. Set
\begin{equation}\label{eq:nonsmooth-hard-oracle}
 \zeta_k=\begin{cases}J_kh_k/\Gamma_k,&\Gamma_k>0,\\0,&\Gamma_k=0,
 \end{cases}
 \qquad
 \widehat{\nabla}_k F(X_k)=\bigl(\gamma_{\mathrm{n}}\sgn(u_k)+\zeta_k,
                 \gamma_{\mathrm{n}}\xi_k\varphi'(v_k)\bigr).
\end{equation}
All $h_k$ are positive. When $\Gamma_k=0$, $\pi_k=0$, so the
zero-noise convention is consistent. One can implement the protocol by
fresh independent uniform random variables each round, applying the
predictable thresholds in \eqref{eq:nonsmooth-protocol-scales} to them.
Thus the construction is measurable and nonanticipating even for a
fully randomized history-dependent policy.

\noindent The thresholds are now explicit functions of pre-sample information. Before
using their probabilities to force recurrent shocks, we check unbiasedness,
the uniform conditional second-moment budget, and the absence of exact
finite-time termination in the second coordinate.

\begin{lemma}\label{lem:nonsmooth-oracle}
For every policy in $\cS_{\mathrm{fin}}$, the oracle satisfies
\begin{equation}\label{eq:nonsmooth-oracle-moments}
 \bbE_k[\widehat{\nabla}_k F(X_k)]\in\partial F(X_k),\qquad \bbE_k[\|{\widehat{\nabla}_k F(X_k)}\|^2]\le49/48<2.
\end{equation}
It can be realized by convex random sample losses with conditional mean
$F$. Moreover,
\begin{equation}\label{eq:nonsmooth-nohit}
 \bbP(v_k\ne0\text{ for every }k\ge1)=1.
\end{equation}
\end{lemma}
\begin{proof}
Symmetry gives $\bbE_k[\zeta_k]=0$, while $\bbE_k[\xi_k]=1$ and
$\bbE_k[\xi_k^2]=13/12$. For positive $\Gamma_k$,
\[
 \bbE_k[\zeta_k^2]=\pi_kh_k^2/\Gamma_k^2\le1/2;
\]
this also holds when $\Gamma_k=0$. The conditional mean is
$(\gamma_{\mathrm{n}}\sgn(u_k),\gamma_{\mathrm{n}}\varphi'(v_k))\in\partial F(X_k)$.
The cross term with $\zeta_k$ has zero conditional expectation, and
$|\varphi'|\le1$. Hence
\[
 \bbE_k[\|{\widehat{\nabla}_k F(X_k)}\|^2]\le\frac14+\frac12+\frac14 \cdot\frac{13}{12}=\frac{49}{48}.
\]
In particular all responses are integrable. The random convex loss
$\gamma_{\mathrm{n}}|u|+\zeta_ku+\gamma_{\mathrm{n}}\xi_k\varphi(v)$ has
conditional expectation $F(u,v)$, and the oracle is its chosen
subgradient at $X_k$. These losses are not asserted to be i.i.d.

\noindent If $v_k\ne0$ and $\Gamma_k=0$, then $v_{k+1}=v_k\ne0$. Otherwise,
projection onto $[-2,2]$ equals zero only when its input is zero, so a
zero next coordinate requires
\[
 v_k-\gamma_{\mathrm{n}}\Gamma_k\xi_k\varphi'(v_k)=0.
\]
Since $\varphi'(v_k)\ne0$, this fixes at most one value of $\xi_k$,
which has probability zero conditionally on $\cF_k$. Starting at $v_1=2$,
induction and a countable union prove \eqref{eq:nonsmooth-nohit}.
\end{proof}

\noindent Oracle admissibility alone does not guarantee a useful shock. Two elementary
facts supply the missing link: the online divisor preserves divergent total
probability, and opposite projected displacements cannot both be small.
Together they handle truncation without requiring it to become inactive.

\begin{lemma}\label{lem:nonsmooth-divisor}
If $a_k\ge0$ and $\sum_k a_k=\infty$, then for
$b_k=\min\{1/2,a_k\}$ and $d_k=\sqrt{1+\sum_{i\le k}b_i}$,
\begin{equation}\label{eq:nonsmooth-divisor}
 d_k\longrightarrow\infty,\qquad\sum_k b_k/d_k=\infty.
\end{equation}
Also, for any $z\in\bbR$ and $0<h\le1$,
\begin{equation}\label{eq:projection-separation}
 \max\{|\proj_{[-1,1]}(z-h)|,|\proj_{[-1,1]}(z+h)|\}\ge h.
\end{equation}
\end{lemma}
\begin{proof}
If $a_k\ge1/2$ infinitely often, then $\sum b_k=\infty$; otherwise $b_k=a_k$ eventually, and the same conclusion holds. Since $d_k^2-d_{k-1}^2=b_k$ for $k\ge2$ and $d_1^2-1=b_1$,
\[
 \frac{b_k}{d_k}\ge d_k-d_{k-1}\quad(k\ge2),\qquad \frac{b_1}{d_1}\ge d_1-1,
\]
so $\sum_{k\le K}b_k/d_k\ge d_K-1\to\infty$ when $K\to\infty$, proving \eqref{eq:nonsmooth-divisor}.
For \eqref{eq:projection-separation}, if both projected magnitudes were
strictly below $h\le1$, neither input could lie outside $[-1,1]$.
Projection would be inactive, so both $z-h$ and $z+h$ would lie in
$(-h,h)$. The first would force $z>0$ and the second $z<0$, a
contradiction.
\end{proof}

\noindent The divisor and projection estimates serve the divergent regime. In the convergent regime, the calibration of the oracle simplifies the slow-movement step of the smooth constructions. There, finite total shock hazard had to be converted into a bound on the cumulative stepsize through a tail estimate and H\"older's inequality (Lemmas~\ref{lem:boundary-finite-hazard} and~\ref{lem:finite-hazard}). Here, the shock hazard is governed directly by the weighted square sum $\sum_kk^{1/2}\Gamma_k^2$, so the convergent regime plays the role of the finite-hazard branch, and the following elementary estimate yields the same conclusion: the cumulative stepsize is little-$o$ of the critical one-quarter scale.

\begin{lemma}\label{lem:weighted-square}
Consider the nonnegative stepsize sequence $(\Gamma_k)_{k\ge1}$ and put
$A_N=\sum_{k=1}^{N-1}\Gamma_k$. Then
\begin{equation}\label{eq:weighted-square-implication}
 \sum_{k\ge1}k^{1/2}\Gamma_k^2<\infty
 \quad\Longrightarrow\quad
 \Gamma_k\to0,\qquad A_N=o(N^{1/4}).
\end{equation}
\end{lemma}
\begin{proof}
The first assertion follows because $k^{1/2}\Gamma_k^2\to0$.
For fixed $m<N$, Cauchy--Schwarz gives
\[
 A_N-A_m
 \le\left(\sum_{k=m}^{N-1}k^{1/2}\Gamma_k^2\right)^{1/2}
     \left(\sum_{k=m}^{N-1}k^{-1/2}\right)^{1/2}
 \le\sqrt2 N^{1/4}
     \left(\sum_{k=m}^{\infty}k^{1/2}\Gamma_k^2\right)^{1/2}.
\]
Divide by $N^{1/4}$, let $N\to\infty$, and then let $m\to\infty$.
The argument is purely numerical, so it applies to every realized path
satisfying the hypothesis.
\end{proof}


\noindent Having verified the oracle and the auxiliary estimates, we prove
the lower bound by separating the two possible weighted-square-sum
regimes. In the divergent regime, the divisor and projection
estimates force recurrent shocks. In the convergent regime,
Lemma~\ref{lem:weighted-square} and the deterministic clock leave a
sufficiently large gap at the sparse times.

\begin{proof}[Theorem~\ref{thm:nonsmooth_lb}]
Fix any policy in $\cS_{\mathrm{fin}}$ and the associated protocol. Its
geometric and oracle properties have just been verified. Define
\[
 \Omega_\infty=\left\{\sum_k k^{1/2}\Gamma_k^2=\infty\right\},\qquad
 \Omega_f=\left\{\sum_k k^{1/2}\Gamma_k^2<\infty\right\}.
\]
They form a measurable partition of $\Omega$.

\paragraph{Divergent weighted square sum.}
On $\Omega_\infty$, Lemma~\ref{lem:nonsmooth-divisor} with $a_k=a_k^{\mathrm n}$ implies $d_k\to\infty$, hence $c_k\to\infty$, and $\sum b_k/d_k=\infty$.
The truncations cause no difficulty: directly from
\eqref{eq:nonsmooth-protocol-scales}, $h_k^2\le d_kk^{-1/2}$, so
\begin{equation}\label{eq:nonsmooth-hazard-lower}
 \pi_k\ge\min\{1/4,a_k^{\mathrm{n}}/(2d_k)\}\ge b_k/(2d_k).
\end{equation}
The last inequality follows from $b_k\le a_k^{\mathrm{n}}$, $b_k\le1/2$,
and $d_k\ge1$. Hence $\sum_k\pi_k=\infty$ on $\Omega_\infty$.
The two nonzero shock signs produce first-coordinate updates
\[
 \proj_{[-1,1]}(\bar u_k-h_k),\qquad
 \proj_{[-1,1]}(\bar u_k+h_k),\qquad
 \bar u_k=u_k-\gamma_{\mathrm{n}}\Gamma_k\sgn(u_k).
\]
By \eqref{eq:projection-separation}, at least one has magnitude at least
$h_k$. Each sign has conditional probability $\pi_k/2$. Therefore
\[
 \bbP_k(|u_{k+1}|\ge h_k)\ge\pi_k/2.
\]
Lemma~\ref{lem:cbc} gives $|u_{k+1}|\ge h_k$ infinitely often, almost
surely on $\Omega_\infty$. At those times,
\[
 (k+1)^{1/4}\Delta(X_{k+1})\ge\gamma_{\mathrm{n}}k^{1/4}h_k
  =\gamma_{\mathrm{n}}\min\{k^{1/4}/2,c_k\}\longrightarrow\infty, \ \ \text{as} \ \ k\to\infty.
\]
This branch already has the infinite limsup without logarithmic
amplification, so it implies the assertion for all $p>0$ at once.

\paragraph{Finite weighted square sum.}
Work on the probability-one no-hit event in
\eqref{eq:nonsmooth-nohit} and fix a path in $\Omega_f$.
Lemma~\ref{lem:weighted-square} gives $\Gamma_k\to0$ and
$A_N=o(N^{1/4})$. Choose a finite path-dependent index $\tau$ such that
$\gamma_{\mathrm{n}}\Gamma_k\xi_k\le1/2$ for every $k\ge\tau$;
this is possible since $\xi_k\le3/2$. Let $w_k=|v_k|>0$.
As $\psi(w)\le w$, these subsequent unprojected steps preserve sign and
have magnitude at most $w_k\le2$. Thus projection is inactive and
\[
 w_{k+1}=w_k-\gamma_{\mathrm{n}}\Gamma_k\xi_k\psi(w_k)
            \in[w_k/2,w_k].
\]
Lemma~\ref{lem:clock} gives
\[
 H_2(w_{k+1})-H_2(w_k)
 \le2\gamma_{\mathrm{n}}\Gamma_k\xi_k\le\frac32\Gamma_k.
\]
Consequently,
\begin{equation}\label{eq:nonsmooth-clock-bound}
 H_2(w_N)\le H_2(w_\tau)+\frac32 A_N=o(N^{1/4}).
\end{equation}
The initial clock value is finite because $w_\tau>0$.

\noindent Use $N=T_j^{\mathrm{n}}$ and $(T_j^{\mathrm{n}})^{1/4}/R_j\to1$, by~\eqref{eq:nonsmooth-scale-growth}.
Equation~\eqref{eq:nonsmooth-clock-bound} and
$H_2(\varrho_j)\ge R_j$ (Lemma~\ref{lem:sparse-profile}) imply $w_{T_j^{\mathrm{n}}}>\varrho_j$ for all
sufficiently large $j$. The objective, the plateau bound of Lemma~\ref{lem:sparse-profile}, and $(T_j^{\mathrm{n}})^{1/4}\ge R_j$ yield
\begin{equation}\label{eq:nonsmooth-subsequence-gap}
 (T_j^{\mathrm{n}})^{1/4}\Delta(X_{T_j^{\mathrm{n}}})
 \ge\frac14(T_j^{\mathrm{n}})^{1/4}\varrho_j\beta_j
 \ge\frac{\varrho_j^2}{4}.
\end{equation}
For every $p>0$, \eqref{eq:nonsmooth-scale-growth} now gives
$[\log(e+T_j^{\mathrm{n}})]^{p}(T_j^{\mathrm{n}})^{1/4}\Delta(X_{T_j^{\mathrm{n}}})\to\infty$, as $j\to\infty$.
Since (77) does not depend on $p$, this holds for all $p>0$ on the same probability-one event.

\noindent The two regimes exhaust every realized weighted-square-sum sequence.
Their combination therefore proves the same infinite-limsup assertion for
the fixed objective and protocol, as required by
Theorem~\ref{thm:nonsmooth_lb}.
\end{proof}

\end{document}